\documentclass[number,final,3p]{elsarticle}
\usepackage{color}
\usepackage{graphicx}
\usepackage{subcaption}
\usepackage{algorithm}
\usepackage[colorlinks]{hyperref}% Add hyper-ref for equations, figures, etc.
\usepackage{amssymb}
\usepackage{amsthm}
\usepackage{amsmath}
\usepackage{booktabs}
\usepackage{url}
\usepackage{scrextend}
\usepackage{epstopdf}
\usepackage{float}
\usepackage{tcolorbox}
\usepackage[perpage]{footmisc}
\usepackage{enumitem}
\usepackage{mathtools}
\usepackage{lineno}
\biboptions{numbers,sort&compress,square}

\journal{arXiv}


\makeatletter
\gdef\urlauthor#1#2{\g@addto@macro\@elsuads{\let\corref\@gobble%
     \def\@@tmp{#1}\raggedright\eadsep
     {\ttfamily\url{\expandafter\strip@prefix\meaning\@@tmp}}\space(#2)%
     \def\eadsep{\unskip,\space}}%
}
\gdef\emailauthor#1#2{\stepcounter{ead}%
     \g@addto@macro\@elseads{\raggedright%
      \let\corref\@gobble\def\@@tmp{#1}%
      \eadsep{\ttfamily\href{mailto:\expandafter\strip@prefix\meaning\@@tmp}{\expandafter\strip@prefix\meaning\@@tmp}}
      (#2)\def\eadsep{\unskip,\space}}%
}
\makeatother

\def\r{\mathbb{R}}
\def\rn{\mathbb{R}^n}
\def\defi{\coloneqq}

\newcommand{\Prob}[1]{\mathbb{P}\left(#1\right)}
\newcommand{\E}[1]{\mathbb{E}\left[#1\right]}
\newcommand{\var}[1]{\mathbb{V}ar\left[#1\right]}
\newcommand{\cov}[2]{\mathbb{C}ov\left[#1\,,#2\right]}
\newcommand{\1}[1]{\mathbb{I}\left(#1\right)}
\newcommand{\vect}[1]{\boldsymbol{#1}}

\begin{document}
\begin{frontmatter}
\renewcommand{\thefootnote}{\alph{footnote}}
\title{Optimized Equivalent Linearization for Random Vibration}
  \author[1]{Ziqi Wang\corref{cor1}}
         \ead{ziqiwang@berkeley.edu}
         \cortext[cor1]{Corresponding author}
  \address[1]{Department of Civil and Environmental Engineering, University of California, Berkeley, United States}
  
\begin{abstract}
A fundamental limitation of various Equivalent Linearization Methods (ELMs) in nonlinear random vibration analysis is that they are approximate by their nature. A quantity of interest estimated from an ELM has no guarantee to be the same as the solution of the original nonlinear system. In this study, we tackle this fundamental limitation. We sequentially address the following two questions: (i) given an equivalent linear system obtained from any ELM, how to construct an estimator such that, as the linear system simulations are guided by a limited number of nonlinear system simulations, the estimator converges on the nonlinear system solution? (ii) how to construct an optimized equivalent linear system such that the estimator approaches the nonlinear system solution as quickly as possible? The first question is theoretically straightforward since classic Monte Carlo techniques, such as the control variates and importance sampling, can improve upon the solution of any surrogate model. We adapt the well-known Monte Carlo theories into the specific context of equivalent linearization. The second question is challenging, especially when rare event probabilities are of interest. We develop specialized methods to construct and optimize linear systems. In the context of uncertainty quantification (UQ), the proposed optimized ELM can be viewed as a \textit{physical surrogate model}-based UQ method. The embedded physical equations endow the surrogate model with the capability to handle high-dimensional uncertainties in stochastic dynamics analysis. 
\end{abstract}


\end{frontmatter}
%\linenumbers

\renewcommand{\thefootnote}{\fnsymbol{footnote}}
\section{Introduction}
\noindent Stochastic dynamic analysis of structural and mechanical systems has been the focus of numerous researchers \cite{RevComptStochMech}\cite{soong_grigoriu_1993}\cite{roberts_spanos_2003}\cite{lutes2004random}\cite{li_chen_2009} in the past several decades. This problem is fundamental in the risk, reliability, and resilience assessment of civil engineering structures subjected to earthquake, hurricane, flood, and tsunami loadings. Linear random vibration problems are relatively straightforward, yet the stochastic dynamic analysis of general multi-degree-of-freedom nonlinear systems still poses significant challenges. Among various nonlinear random vibration methods, the equivalent linearization method (ELM) is a classical approach with a history of 70 years. The topic finds its root in the early works of Booton (1953) \cite{booton1953analysis} and Caughey (1963) \cite{caughey1963equivalent}, and later it gained wide popularity, particularly among engineers. This popularity is partly because the method is easy to implement, efficient, and applicable to relatively general nonlinear systems. Moreover, the equivalent linear system established in an ELM has clear engineering/physical interpretations. Other classical random vibration methods, Fokker-Planck equation, stochastic averaging, moment closure, and perturbation (see \cite{lutes2004random} for a review), although can be more accurate, are often restricted to specialized nonlinear systems. Recent progress on these classical random vibration methods, especially on stochastic differential equation-based methods \cite{chen2017beating}\cite{petromichelakis2020addressing}\cite{chen2016gf}, have made their applications to large-scale systems possible. However, they are still subjected to limitations such as high computational demand, complexity in implementation, weak scalability toward higher-dimensional problems, and lack of interpretability. It is futile to argue which nonlinear random vibration method is the ``optimal"; there is no free lunch: a trade-off between accuracy, efficiency, and generality is inevitable. This work focuses on further developing ELM, aiming to overcome its limitation at an acceptable computational cost.

In ELM, the nonlinear system is replaced by an ``equivalent" linear system. Consequently, the well-developed linear system analysis and linear random vibration theories can be applied to estimate various response quantities. In the conventional ELM \cite{roberts_spanos_2003}\cite{crandall2006half}\cite{elishakoff2017sixty}\cite{atalik1976stochastic}, parameters of the equivalent linear system are determined via minimizing the mean-square error between the responses of the nonlinear and linear systems. The limitation is that the estimated response distributions and rare event probabilities can be far from correct. In a recent development of ELM, namely the Tail-Equivalent Linearization Method (TELM) \cite{fujimura2007tail}\cite{broccardo2014further}, a nonparametric linear system is constructed through a discretized impulse/frequency response function, using knowledge of the ``design point" obtained from the first-order reliability method (FORM). As the name suggests, TELM has good accuracy in estimating the response distributions in the tail region at the cost of performing FORM analysis. In another recent development, the Gaussian Mixture-based ELM (GM-ELM) \cite{wang2017equivalent}\cite{yi2018bivariate}\cite{yi2019gaussian}, the nonlinear system response is decomposed into a probabilistic superposition of linear system responses. The GM-ELM can easily capture the non-Gaussianity of nonlinear system responses due to the flexibility of Gaussian mixture models. However, GM-ELM lacks a clear (classical) physical interpretation, making it not straightforward to integrate with linear system analysis and random vibration techniques.

Despite the recent progress in ELM, the fundamental limitation\textemdash equivalent linearization is approximate by its very nature\textemdash remains untouched. Unless the nonlinearity is negligible, a quantity of interest estimated from an ELM has no guarantee to be the same as the solution of the nonlinear system. In this work, we tackle this fundamental limitation of ELM by sequentially addressing the following two questions.
\begin{itemize}
\item Given an equivalent linear system obtained from any ELM, how to construct an estimator such that, as the linear system simulations are guided by a limited number of nonlinear system simulations, the estimator converges on the nonlinear system solution?
\item How to construct an optimized equivalent linear system such that the estimator approaches the nonlinear system solution as quickly as possible?
\end{itemize}

The well-known Monte Carlo techniques such as the control variates and importance sampling (see \cite{rubinstein2016simulation} for a comprehensive review) can be leveraged to answer the first question. We adapt the Monte Carlo theories into the specific context of ELM, with the technical details introduced in Section \ref{Sec:2}. The second question is challenging and requires some specialized techniques. We address it in Section \ref{Sec:3} and Section \ref{Sec:4}. Conceptually, the main procedures involve: (i) optimizing a parametric linear system such that within some critical region,  the linear system response is maximally correlated with the nonlinear system response; (ii) coupling the linear system simulation with the nonlinear system simulation so that the linear system solution can be corrected toward the nonlinear system solution. 

The general principle that underpins the proposed ELM is to combine simplified/analytic physical models with numerical gap-filling schemes to estimate quantities of complex physical models. In structural reliability, the model correction factor method \cite{ModelCorrect}\cite{xu2015reliability} is one of the earliest approaches that embody the aforementioned principle. In stochastic dynamics, variants of the probability density evolution method \cite{ZHOU2022109435}\cite{ZHOU2022108283} have been developed to combine analytic approximations with numerical schemes. The closest analogy of the proposed method can be found in the field of surrogate modeling-based uncertainty quantification (UQ) \cite{sudret2012meta}\citep{dhulipala2022active}. Specifically, the optimized ELM can be viewed as a \textit{physical surrogate model}-based UQ method. Classical surrogate modeling techniques such as the Gaussian process/Kriging \cite{williams1998prediction} and polynomial chaos expansion \cite{xiu2010numerical} would encounter significant difficulties in solving stochastic dynamics problems, because random processes or their discretized representations are infinite/high-dimensional. Deep neural network \cite{lecun2015deep} and physics-informed neural network \citep{raissi2019physics} can handle high-dimensional input-output of dynamic models. However, their training (a highly complex non-convex optimization) and interpretability still pose open challenges. On the other hand, a physics-based surrogate model is interpretable by construction. The embedded physical equations endow the surrogate model with the capability to handle high-dimensional uncertainties in stochastic dynamics analysis. The idea of leveraging a cheap physics-based model in the probabilistic analysis is also echoed in the multi-fidelity UQ framework \cite{peherstorfer2018survey}. Using the terminology of multi-fidelity UQ, the optimized linear system is a ``low-fidelity model", and we develop specialized ``model fusion and filtering" techniques to couple the linear system simulations with a limited number of nonlinear system simulations. Therefore, the proposed method can also be classified as a specialized multi-fidelity UQ approach for nonlinear stochastic dynamics analysis. Finally, despite the possible connections to various probabilistic analysis methods, our initial and primary motivation remains the exploration of new directions for the equivalent linearization methodology. 

\section{Pushing Equivalent Linear System Solutions toward the Nonlinear System Solutions}\label{Sec:2}
\subsection{Quantities of Interest and Equivalent Linear Systems}\label{Sec:QoI}
\noindent Let a vector $\vect X=[X_1,...,X_n]$ represent the basic\footnote{Here ``basic" indicates that the random variables represent the source of randomness. The source of randomness can be a vague idea; it becomes unambiguous when a model universe is clearly specified \cite{der2009aleatory}. The source of randomness in random vibration analysis is typically stochastic excitation.} random variables considered in a random vibration problem, and let $f_{\vect X}(\vect x)$ denote the joint probability density function of $\vect X$. We assume the stochastic excitation process has been discretized by a finite set of random variables. Consequently, hereafter we will work on the finite-dimensional probability space $(\rn,\mathcal{B}_n,\mathbb{P}_{\vect X})$, where $\mathcal{B}_n$ is the Borel $\sigma$-algebra on $\rn$, and $\mathbb{P}_{\vect X}$ is the probability measure of $\vect X$ with the distribution function $f_{\vect X}(\vect x)$. If the stochastic excitation is a wide-band process, a large number (e.g., $>100$) of random variables is required for an accurate discretization.  

Now consider a computational model, $\mathcal{M}$, that maps the high-dimensional $\vect X$ into a Quantity of Interest (QoI) $Q$, i.e., $\mathcal{M}:\vect X\in\rn\mapsto Q\in\r$. In nonlinear random vibration analysis, $\mathcal{M}$ is an implicit function of $\vect X$ and involves solving nonlinear dynamic equations. The QoI $Q$ is also a random variable, and in practice, it is often of interest to estimate the expectation $\E{Q}$. It is useful to observe that a probability of interest can also be cast into the form of $\E{Q}$. For a generic probability estimate, one could define $\mathcal{M}$ as $\mathcal{M}(\vect X)=\1{\vect X\in\Omega}$, where $\mathbb{I}$ is a binary indicator function and $\Omega$ is a measurable set in $\rn$; consequently, $\E{Q}=\mathbb{P}_{\vect X}{\left(\Omega\right)}$. Given that $\E{Q}$ exists, the integral to estimate $\E{Q}$ is
\begin{equation}\label{IntQ}
\E{Q}=\int_{\vect x\in\rn}\mathcal{M}(\vect x)f_{\vect X}(\vect x)\,d\vect x\,.
\end{equation}
This integral can be fairly challenging because $\vect x$ is high-dimensional, and $\mathcal{M}$ involves solving nonlinear dynamic equations.   

Let $Q_L=\mathcal{M}_L(\vect X;\vect\theta)$ represent the QoI mapping of an equivalent linear system with parameter set $\vect\theta$. Typically, $\vect\theta$ is associated with the mass, stiffness, and damping matrices of the linear system. It is important to note that in general, $\mathcal{M}_L$ is not a linear function of $\vect X$: a linear system can admit nonlinear functions for specific $(\vect X, Q_L)$ pairs. For a concrete example, let $Q_L$ be the peak of a specified response quantity $Y(t)$ of a linear system over some time period $T$, i.e.,  $Q_L=\sup_{t\in T}Y(t)$, and let $\vect X$ correspond to a stochastic excitation $Z(t)$ with the form $Z(t)=\vect s(t)\cdot\vect X$, where $\vect s(t)$ are basis functions depending on the specific random process discretization method \cite{shinozuka1972digital}\cite{chatfield2013analysis}; it follows that $\mathcal{M}_L(\vect X)=\sup_{t\in T}\left[\mathcal{I}(t)\ast(\vect s(t)\cdot\vect X)\right]$, where $\mathcal{I}(t)$ is the impulse response function for the $(Z(t),Y(t))$ pair, and ``$\ast$" denotes convolution. 

The expected QoI of the equivalent linear system, $\E{Q_L}$, is estimated from
\begin{equation}\label{IntQL}
\E{Q_L}=\int_{\vect x\in\rn}\mathcal{M}_L(\vect x;\vect\theta)f_{\vect X}(\vect x)\,d\vect x\,.
\end{equation}
Ideally, $\E{Q_L}$ should approximate $\E{Q}$. Since linear systems are cheap to simulate, we consider $\E{Q_L}$ as known when an equivalent linear system is specified. In the subsequent discussions, the parameter set $\vect\theta$ will be explicitly shown only when necessary.

Given an equivalent linear system specified, and $\mu_L\defi\E{Q_L}$ known, subsequently, we will adopt and adapt two classic Monte Carlo techniques to push the equivalent linear system solution toward the nonlinear system solution. 

\subsection{The Method of Control Variates}\label{Sec:CV}
\noindent Using the method of control variates \cite{rubinstein2016simulation}, we construct the following random variable
\begin{equation}\label{ConVar}
Q_C=Q-\alpha(Q_L-\mu_L)\,,
\end{equation}
where $\alpha\neq0$ is a parameter. Clearly, $\E{Q_C}=\E{Q}$, and the variance of $Q_C$ is
\begin{equation}\label{VarConVar}
\var{Q_C}=\var{Q}+\alpha^2\var{Q_L}-2\alpha\cov{Q}{Q_L}\,,
\end{equation}
where $\mathbb{C}ov$ denotes the covariance. An optimal parameter setting of $\alpha$ can be found by setting the derivative $\frac{d\,\var{Q_C}}{d\,\alpha}$ to $0$, and the optimal $\alpha$ is
\begin{equation}\label{alpha}
\alpha^*=\dfrac{\cov{Q}{Q_L}}{\var{Q_L}}\,,
\end{equation}
where the variance of $Q_L$ is also considered as known and is denoted by $\sigma_L^2\defi\var{Q_L}$ for later use. Substituting Eq.\eqref{alpha} back into Eq.\eqref{VarConVar}, we get
\begin{equation}\label{OptVarConVar}
\mathbb{V}ar^*\left[Q_C\right]=(1-\rho^2)\var{Q}\,,
\end{equation}
where the correlation coefficient $\rho$ is
\begin{equation}\label{Correl}
\rho=\dfrac{\cov{Q}{Q_L}}{\sqrt{\var{Q}\var{Q_L}}}\,.
\end{equation}
Eq.\eqref{OptVarConVar} suggests that if $\E{Q_C}$ is used to estimate $\E{Q}$, there will be a variance reduction as long as $Q$ and $Q_L$ are correlated (negatively or positively). Moreover, the variance reduction is irrelevant to the magnitude of $Q_L$, i.e., introducing a scale factor to $Q_L$ will leave the variance reduction intact. As long as $Q_L$ is highly correlated with $Q$, i.e., $|\rho|\approx1$, the method would be highly efficient since $\var{Q_C}\approx0$. This observation instantly suggests that the equivalent linear system can be optimized to yield $|\rho|\approx1$ and thus maximum efficiency. This idea will be explored later. 

In practice, the optimal $\alpha^*$ cannot be obtained since $\cov{Q}{Q_L}$ is unknown. However, it is straightforward to build a random sampling-based numerical method to approximate $\alpha^*$ and consequently $\E{Q_C}$. The flowchart of the control variate-enhanced ELM is { developed in \ref{Fig:CV}}. Note that to improve the readability of this paper, we list all the implementation details in the appendix, while the main text focuses on the methodological development and its insight. 

The control variate-enhanced ELM is straightforward to implement; this is particularly attractive to engineers. However, the limitation is that the approximation relies on sampling from $f_{\vect X}(\vect x)$. If most samples drawn from $f_{\vect X}(\vect x)$ do not contribute much to $\E{Q}$, the method would be ineffective. {Moreover, the correlation coefficient $\rho$ can be a crude metric for extreme responses; consequently, the control variate-enhanced ELM may only achieve a limited efficiency gain for extreme value estimations. These limitations motivate the study of the following method.}

\subsection{The Method of Linear System-based Importance Sampling}\label{Sec:OIS}
\noindent Let an importance density, $h(\vect x)$, be introduced to estimate Eq.\eqref{IntQ}, and the integral is rewritten as
\begin{equation}\label{IS}
\E{Q}=\int_{\vect x\in\rn}\mathcal{M}(\vect x)\dfrac{f_{\vect X}(\vect x)}{h(\vect x)}h(\vect x)\,d\vect x=\mathbb{E}_h\left[\mathcal{M}(\vect X)\dfrac{f_{\vect X}(\vect X)}{h(\vect X)}\right]\,.
\end{equation}
The subscript of $\mathbb{E}_h$ highlights the expectation is taken with respect to the importance density $h(\vect x)$ instead of $f_{\vect X}(\vect x)$. In the importance sampling practice, if the dimensionality of $\vect x$ is large, which is the case of random vibration, constructing $h(\vect x)$ is fairly challenging. Despite the practical challenge, there exists a theoretical optimal importance density \cite{rubinstein2016simulation}, written as
\begin{equation}\label{PerfectIS}
h^*(\vect x)=\dfrac{\vert\mathcal{M}(\vect x)\vert f_{\vect X}(\vect x)}{\int_{\vect x\in\rn}\vert\mathcal{M}(\vect x)\vert f_{\vect X}(\vect x)\,d\vect x}\,.
\end{equation}
In particular, if $\forall\vect x$, $\mathcal{M}(\vect x)\geq0$ or $\mathcal{M}(\vect x)\leq0$, the optimal importance density leads to a zero variance estimation. The catch of the optimal density is that the estimation of the normalizing constant is as challenging as (sometimes\footnote{By sometimes, we mean $\forall\vect x$, $\mathcal{M}(\vect x)\geq0$ or $\mathcal{M}(\vect x)\leq0$.} exactly equivalent to) the estimation of $\E{Q}$.

Replacing the $\mathcal{M}(\vect x)$ in Eq.\eqref{PerfectIS} by $\mathcal{M}_L(\vect x)$ of the equivalent linear system, we obtain
\begin{equation}\label{PerfectISLin}
h_L^*(\vect x)=\dfrac{\vert\mathcal{M}_L(\vect x)\vert f_{\vect X}(\vect x)}{\int_{\vect x\in\rn}\vert\mathcal{M}_L(\vect x)\vert f_{\vect X}(\vect x)\,d\vect x}=\dfrac{\vert\mathcal{M}_L(\vect x)\vert f_{\vect X}(\vect x)}{\mu_{|L|}}\,.
\end{equation}
Clearly, the density $h_L^*(\vect x)$ is optimal to estimate $\E{Q_L}$, i.e., zero-variance to estimate $\E{Q_L}$. The normalizing constant $\mu_{|L|}$ only involves the linear system and thus is considered as known. In the special case of $\forall\vect x$, $\mathcal{M}(\vect x)\geq0$ or $\mathcal{M}(\vect x)\leq0$, we have $\mu_{|L|}=\mu_{L}$ or $\mu_{|L|}=-\mu_{L}$, respectively.

Replacing the $h(\vect x)$ in Eq.\eqref{IS} by $h_L^*(\vect x)$, we obtain
\begin{equation}\label{ISLin}
\E{Q}=\mathbb{E}_{h_L^*}\left[\mu_{|L|}\dfrac{\mathcal{M}(\vect X)}{\vert\mathcal{M}_L(\vect X)\vert}\right]\,.
\end{equation}
The term $\frac{\mathcal{M}(\vect X)}{\vert\mathcal{M}_L(\vect X)\vert}$ can be interpreted as a correction factor for the linear system solution $\mu_{|L|}$, and the expectation $\mathbb{E}_{h_L^*}\left[\frac{\mathcal{M}(\vect X)}{\vert\mathcal{M}_L(\vect X)\vert}\right]$ serves as an averaged correction over $\mu_{|L|}$. 

The underlying assumption of using Eq.\eqref{ISLin} is: if the equivalent linear system resembles the nonlinear system, the optimal importance density for $\E{Q_L}$ can be an effective importance density for $\E{Q}$. However, it is important to note that even if the  linear system response does not resemble the nonlinear response, the importance sampling theory still guarantees the correctness of Eq.\eqref{ISLin}. Specifically, as long as $h_L^*(\vect x)$ \textit{dominates} $\mathcal{M}(\vect x)f_{\vect X}(\vect x)$, i.e., $h_L^*(\vect x)=0\implies\mathcal{M}(\vect x)f_{\vect X}(\vect x)=0$, the importance sampling via $h_L^*(\vect x)$ provides an \textit{unbiased} estimator of $\E{Q}$. The dominance condition of $h_L^*(\vect x)$ can be met if $\mathcal{M}_L(\vect x)$ is non-zero \textit{almost everywhere}\footnote{For the region where $\mathcal{M}(\vect x)f_{\vect X}(\vect x)=0$, the value of $\mathcal{M}_L(\vect x)$ can be zero (or any real number).}. On the other hand, if $\mathcal{M}(\vect x)$ and $\mathcal{M}_L(\vect x)$ are binary indicator functions, i.e., $\E{Q}$ represents a probability, specialized adaptations are required. We will address this issue in a separate section. Despite the theoretical correctness, the efficiency of Eq.\eqref{ISLin} does rely on the ``closeness" between $h_L^*(\vect x)$ and $h^*(\vect x)$; this motivates the optimization of $h_L^*(\vect x)$ that will be studied in the following section.

A major technical challenge of using Eq.\eqref{ISLin} lies in the sampling of $h_L^*(\vect x)$. This challenge can be handled by various Markov Chain Monte Carlo (MCMC) samplers \cite{andrieu2003introduction}. Note that sampling from $h_L^*(\vect x)$ via an MCMC sampler only involves simulations of the equivalent linear system.

Finally, it is useful to derive the following variance. 
\begin{equation}\label{VarIS}
\begin{aligned}
\mathop{\mathbb{V}ar_{h_L^*}}\left[\mu_{|L|}\dfrac{\mathcal{M}(\vect X)}{\vert\mathcal{M}_L(\vect X)\vert}\right]&=\int_{\vect x\in\rn}\left(\mu_{|L|}\dfrac{\mathcal{M}(\vect x)}{\vert\mathcal{M}_L(\vect x)\vert}\right)^2\frac{\vert\mathcal{M}_L(\vect x)\vert f_{\vect X}(\vect x)}{\mu_{|L|}}\,d\vect x-\left(\E{Q}\right)^2\\
&=\int_{\vect x\in\rn}\mathcal{M}^2(\vect x)\dfrac{\mu_{|L|}}{\vert\mathcal{M}_L(\vect x)\vert}f_{\vect X}(\vect x)\,d\vect x-\left(\int_{\vect x\in\rn}\mathcal{M}^2(\vect x)f_{\vect X}(\vect x)\,d\vect x-\var{Q}\right)\\
&=\var{Q}-\E{\mathcal{M}^2(\vect X)\left(1-\dfrac{\mu_{|L|}}{\vert\mathcal{M}_L(\vect X)\vert}\right)}\\
&=\var{Q}-\cov{\dfrac{Q^2}{\vert Q_L\vert}}{\vert Q_L\vert}\,.
\end{aligned}
\end{equation}
Eq.\eqref{VarIS} suggests that the equivalent linear system-based importance sampling formula, Eq.\eqref{ISLin}, encodes a variance reduction as long as $\cov{\frac{Q^2}{\vert Q_L\vert}}{\vert Q_L\vert}>0$, i.e., $\frac{Q^2}{\vert Q_L\vert}$ and $\vert Q_L\vert$ are positively correlated. In the special case where $Q=a|Q_L|$, $a\neq0$, the importance sampling estimator is zero-variance. Similar to the control variate-enhanced ELM, Eq.\eqref{VarIS} also suggests that the magnitude of $Q_L$ is irrelevant to the variance reduction. The implementation details of the importance sampling-enhanced ELM are {developed in \ref{Fig:OIS}}.

\section{Optimized Equivalent Linear Systems}\label{Sec:3}
\subsection{Parametrization}
\noindent We model the equivalent linear system by a generic impulse response function expressed as
\begin{equation}\label{LinearSystem}
\mathcal{I}(t;\vect\theta)=b\sum_{i=1}^{DoF}a_i\,\mathcal{I}_i(t;\omega_i,\zeta_i)\,,
\end{equation}
where $\mathcal{I}_i(t;\omega_i,\zeta_i)$ is
\begin{equation}\label{LinearSystemi}
\mathcal{I}_i(t;\omega_i,\zeta_i)=\dfrac{1}{\omega_i\sqrt{1-\zeta_i^2}}e^{-\zeta_i\omega_it}\sin\left(\omega_i\sqrt{1-\zeta_i^2}\,t\right)\,.
\end{equation}

The full parameter set is $\lbrace\vect\theta\rbrace=\lbrace DoF,b,a_1,\omega_1,\zeta_1,\dots,a_{DoF},\omega_{DoF},\zeta_{DoF}\rbrace$. We summarize the parameters and their meanings in the following table.

\begin{table}[H]
  \caption{\textbf{Parameters of the Equivalent Linear System}}
  \label{tab:ELM}
  \centering
  \begin{tabular}{l l l}
    \toprule
    Parameter & Description & Remark\\
    \midrule
    $DoF$     & Degree-of-Freedom of the linear system & controls the modeling flexibility\\
    $b$       & global scale factor of the linear response & {often $=1$ (more remarks below)}\\
    $a_i$     & weight of the $i$-th modal response & constrained by $\sum_{i}a_i^2=1$ \\
    $\omega_i$ & natural frequency of the $i$-th mode               & \\
    $\zeta_i$ & damping ratio of the $i$-th mode               & \\
    \bottomrule
  \end{tabular}
\end{table}
Additional remarks on the parameters are listed below. 

\begin{itemize}[leftmargin=*]
\item \textbf{The parameter $DoF$}\\
The $DoF$ of the equivalent linear system can be irrelevant to that of the original nonlinear system. In the context of optimization, $DoF$ is treated similarly to other parameters in $\vect\theta$, i.e., $DoF$ is determined through optimization. 
On top of the fact that $DoF$ must be a positive integer, the significant difference between $DoF$ and other parameters is that there is a hierarchical structure: $DoF$ determines the number of $a_i$, $\omega_i$, and $\zeta_i$. Therefore, a double layer optimization can be adopted, i.e., to iteratively increase $DoF$ and solve the optimization conditional on each $DoF$ value. The optimal $DoF$ is found if a performance measure does not significantly increase.
\item \textbf{The parameter $b$}\\
Implied from Section \ref{Sec:CV} and Section \ref{Sec:OIS}, the global scale factor $b$ would not impact the variance reduction and can be fixed to 1, as long as the influence of $b$ within the model $\mathcal{M}_L$ can be factorized out. For a concrete example, recall the peak response model $\mathcal{M}_L(\vect X)=\sup_{t\in T}\left[\mathcal{I}(t)\ast(\vect s(t)\cdot\vect X)\right]$ introduced in Section \ref{Sec:QoI}. If $\mathcal{I}(t)$ is scaled by $b\mathcal{I}(t)$, the parameter $b$ can be factorized out of $\sup(\cdot)$, and $\mathcal{M}_L(\vect X;b)=b\mathcal{M}_L(\vect X)$. On the other hand, if $\mathcal{M}_L(\vect X)=\1{\sup_{t\in T}\left[b\mathcal{I}(t)\ast(\vect s(t)\cdot\vect X)\right]\in\Omega}$, the factor $b$ cannot be factorized out of the indicator function. 
\item \textbf{Parameters $a_i$}\\
The constraint of $a_i$, $\sum_{i}a_i^2=1$, is artificially designed to reduce the parametric redundancy and facilitate optimization. If $a_i$ is not normalized, the effect of $b$ and $a_i$ will overlap.
\item \textbf{Parameters $\omega_i$}\\
To be physically admissible, one should set $\omega_i>0$. In practice, one may consider a cut-off frequency $\omega_{\max}$ for the dynamical analysis, and thus one can set $\omega_i\in(0,\omega_{\max}]$.
\item \textbf{Parameters $\zeta_i$}\\
To be physically admissible, one should set $\zeta_i\geq0$. In this study, we further restrict to the \textit{underdamped} system, i.e., $\zeta_i\in[0,1)$. The extensions to critically damped and overdamped systems are straightforward, given that the impulse response function Eq.\eqref{LinearSystemi} is modified accordingly.
\end{itemize}

\subsection{Optimization for the Control Variate-Enhanced ELM}
\noindent Suggested by Eq.\eqref{OptVarConVar}, the optimization to be solved in the control variate-enhanced ELM is
\begin{equation}\label{OptCV}
\vect\theta^*=\mathop{\arg\max}_{\vect\theta}\rho^2(\vect\theta)\,.
\end{equation}

In reality, $\rho$ can only be approximated using random samples, and the optimization becomes
\begin{equation}\label{OptCVP}
\vect\theta^*\approx\mathop{\arg\max}_{\vect\theta}\hat\rho^2(\vect\theta|\mathcal{\vect D})\,,
\end{equation}
where the dataset $\mathcal{\vect D}$ contains random samples of $\vect X$ and the corresponding $Q=\mathcal{M}(\vect X)$. Due to the wide spectrum of possible QoIs, we do not expect Eq.\eqref{OptCVP} to be solved analytically. Therefore, in general, Eq.\eqref{OptCVP} should be solved by various derivative-free optimization techniques \cite{conn2009introduction}. Note that for a given $\mathcal{\vect D}$, solving Eq.\eqref{OptCVP} does not introduce additional nonlinear system simulations. The flowchart of the adaptive control variate-enhanced ELM is developed in \ref{Fig:ACV}.

\subsection{Optimization for the Importance Sampling-Enhanced ELM}\label{Sec:OOIS}
\noindent Suggested by Eq.\eqref{VarIS}, the optimization to be solved in the importance sampling-enhanced ELM is
\begin{equation}\label{OptIS}
\begin{aligned}
\vect\theta^*&=\mathop{\arg\min}_{\vect\theta}\left\lbrace\E{\dfrac{Q^2}{\vert Q_L(\vect\theta)\vert}}\E{\vert Q_L(\vect\theta)\vert}\right\rbrace\\
&=\mathop{\arg\min}_{\vect\theta}\left\lbrace\int_{\vect x\in\rn}\dfrac{\mathcal{M}^2(\vect x)}{\vert \mathcal{M}_L(\vect x;\vect\theta)\vert}f_{\vect X}(\vect x)\,d\vect x\,\int_{\vect x\in\rn}\vert \mathcal{M}_L(\vect x;\vect\theta)\vert f_{\vect X}(\vect x)\,d\vect x\right\rbrace\,.
\end{aligned}
\end{equation}
Notice that Eq.\eqref{OptIS} attempts to minimize the sampling variance by maximizing the covariance in Eq.\eqref{VarIS}, dropping the constant component independent of $\vect\theta$. The second integral inside the curled brackets of Eq.\eqref{OptIS}, alternatively denoted by $\mu_{|L|}(\vect\theta)$, is easy to estimate since it only involves the equivalent linear system; while the first integral is challenging since it inevitably requires simulations of the nonlinear system. For computational efficiency, it is useful to introduce importance sampling into Eq.\eqref{OptIS}. Specifically, an iterative approach can be designed to find $\vect\theta^*$, and the importance density can be adaptively set as the previous solution of Eq.\eqref{OptIS}. Formally, let $\vect\theta^*_{j}$ denote the optimal parameter set to be estimated at the current iterative step, and $\vect\theta^*_{j-1}$ denote the known solution from the previous step. The optimization Eq.\eqref{OptIS} is rewritten as
\begin{equation}\label{OptISj}
\vect\theta^*_{j}=\mathop{\arg\min}_{\vect\theta}\left\lbrace\int_{\vect x\in\rn}\dfrac{\mathcal{M}^2(\vect x)}{\vert \mathcal{M}_L(\vect x;\vect\theta)\vert}\dfrac{f_{\vect X}(\vect x)}{h^*_L(\vect x;\vect\theta^*_{j-1})}h^*_L(\vect x;\vect\theta^*_{j-1})\,d\vect x\,\int_{\vect x\in\rn}\dfrac{\vert\mathcal{M}_L(\vect x;\vect\theta)\vert f_{\vect X}(\vect x)}{h^*_L(\vect x;\vect\theta^*_{j-1})}h^*_L(\vect x;\vect\theta^*_{j-1})\,d\vect x\right\rbrace\,,
\end{equation}
where the importance density $h^*_L(\vect x;\vect\theta^*_{j-1})$ is
\begin{equation}\label{ISj}
h^*_L(\vect x;\vect\theta^*_{j-1})=\dfrac{\vert\mathcal{M}_L(\vect x;\vect\theta^*_{j-1})\vert f_{\vect X}(\vect x)}{\int_{\vect x\in\rn}\vert\mathcal{M}_L(\vect x;\vect\theta^*_{j-1})\vert f_{\vect X}(\vect x)\,d\vect x}=\dfrac{\vert\mathcal{M}_L(\vect x;\vect\theta^*_{j-1})\vert f_{\vect X}(\vect x)}{\mu_{|L|}(\vect\theta^*_{j-1})}\,.
\end{equation}
Substituting Eq.\eqref{ISj} back into Eq.\eqref{OptISj}, we obtain
\begin{equation}\label{OptISjF}
\begin{aligned}
\vect\theta^*_{j}&=\mathop{\arg\min}_{\vect\theta}\left\lbrace\int_{\vect x\in\rn}\dfrac{\mathcal{M}^2(\vect x)}{\vert \mathcal{M}_L(\vect x;\vect\theta)\mathcal{M}_L(\vect x;\vect\theta^*_{j-1})\vert}h^*_L(\vect x;\vect\theta^*_{j-1})\,d\vect x\,\int_{\vect x\in\rn}\Big|\dfrac{\mathcal{M}_L(\vect x;\vect\theta)}{\mathcal{M}_L(\vect x;\vect\theta^*_{j-1})}\Big|h^*_L(\vect x;\vect\theta^*_{j-1})\,d\vect x\right\rbrace\\
&\approx\mathop{\arg\min}_{\vect\theta}\left\lbrace\left\langle\frac{\mathcal{M}^2(\vect x)}{\vert\mathcal{M}_L(\vect x;\vect\theta)\mathcal{M}_L(\vect x;\vect\theta^*_{j-1})\vert}\right\rangle_{h^*_L}\left\langle\Big|\dfrac{\mathcal{M}_L(\vect x;\vect\theta)}{\mathcal{M}_L(\vect x;\vect\theta^*_{j-1})}\Big|\right\rangle_{h^*_L}\right\rbrace\,,
\end{aligned}
\end{equation}
where $\langle\cdot\rangle_{h^*_L}$ denotes the sample mean computed from random samples of $h^*_L(\vect x;\vect\theta^*_{j-1})$. Note that $\mu_{|L|}(\vect\theta^*_{j-1})$ is dropped from Eq.\eqref{OptISjF} since $\vect\theta^*_{j-1}$ is fixed at the iterative step $j$. Also, note that the same importance density is introduced to estimate $\mu_{|L|}(\vect\theta)$. This is to enforce the approximations of the two integrals being consistent, and the errors could, to some extent, cancel out.

The motivation of Eq.\eqref{OptIS} is to minimize the variance. An alternative chain of reasoning to construct an optimization is: (i) the optimal importance density is described by Eq.\eqref{PerfectIS}; (ii) an ideal linear system should lead to a density $h^*_L(\vect x;\vect\theta^*)$ being similar to Eq.\eqref{PerfectIS}; (iii) the optimal $\vect\theta^*$ can be found by minimizing a discrepancy measure between $h^*_L(\vect x;\vect\theta^*)$ and Eq.\eqref{PerfectIS}. If the cross entropy is adopted as a discrepancy measure, the optimization has the following form,
\begin{equation}\label{OptCE}
\begin{aligned}
\vect\theta^*&=\mathop{\arg\min}_{\vect\theta}\mathop{CE}\left(h^*(\vect x),h^*_L(\vect x;\vect\theta)\right)\\
&=\mathop{\arg\min}_{\vect\theta}\left\lbrace-\int_{\vect x\in\rn}h^*(\vect x)\log h^*_L(\vect x;\vect\theta)\,d\vect x\right\rbrace\\
&=\mathop{\arg\max}_{\vect\theta}\int_{\vect x\in\rn}|\mathcal{M}(\vect x)|\log\frac{|\mathcal{M}_L(\vect x;\vect\theta)|}{\mu_{|L|}(\vect\theta)}f_{\vect X}(\vect x)\,d\vect x\,.
\end{aligned}
\end{equation}
Similarly, an adaptive importance sampling can be introduced to solve Eq.\eqref{OptCE}. The optimization to be solved in the sequence is written as
\begin{equation}\label{OptCEj}
\begin{aligned}
\vect\theta^*_{j}&=\mathop{\arg\max}_{\vect\theta}\int_{\vect x\in\rn}|\mathcal{M}(\vect x)|\log\frac{|\mathcal{M}_L(\vect x;\vect\theta)|}{\mu_{|L|}(\vect\theta)}\dfrac{f_{\vect X}(\vect x)}{h^*_L(\vect x;\vect\theta^*_{j-1})}h^*_L(\vect x;\vect\theta^*_{j-1})\,d\vect x\\
&=\mathop{\arg\max}_{\vect\theta}\int_{\vect x\in\rn}\Big|\dfrac{\mathcal{M}(\vect x)}{\mathcal{M}_L(\vect x;\vect\theta^*_{j-1})}\Big|\log\frac{|\mathcal{M}_L(\vect x;\vect\theta)|}{\mu_{|L|}(\vect\theta)}h^*_L(\vect x;\vect\theta^*_{j-1})\,d\vect x\\
&\approx\mathop{\arg\max}_{\vect\theta}\left\langle\Big|\dfrac{\mathcal{M}(\vect x)}{\mathcal{M}_L(\vect x;\vect\theta^*_{j-1})}\Big|\Big(\log|\mathcal{M}_L(\vect x;\vect\theta)|-\log\Big\langle\Big|\dfrac{ \mathcal{M}_L(\vect x;\vect\theta)}{\mathcal{M}_L(\vect x;\vect\theta^*_{j-1})}\Big|\Big\rangle_{h^*_L}\Big)\right\rangle_{h^*_L}\,.
\end{aligned}
\end{equation}
Similar to Eq.\eqref{OptISjF}, the importance density is also introduced to estimate $\mu_{|L|}(\vect\theta)$, with the motivation of consistent accuracy. The flowchart of the adaptive importance sampling-enhanced ELM is {developed in \ref{Fig:AOIS}}.

\section{Adaptations for Rare Event Probability Estimation}\label{Sec:4}
\subsection{Linear System-based Importance Density}
\noindent Eq.\eqref{IntQ} is particularly challenging if $\E{Q}$ represents a rare event probability, i.e., $Q=\1{\vect X\in\Omega}$ so that $\E{Q}=\mathbb{P}_{\vect X}(\Omega)$, and the probability is small. This is the typical setting in reliability analysis and catastrophic risk analysis. In this context, the control variate-enhanced ELM will be inefficient because most samples generated from $f_{\vect X}(\vect x)$ would not contribute to the rare event probability. Therefore, this section will focus on the importance sampling-enhanced ELM.

To be consistent with the notations in previous sections, consider $Q=\1{\mathcal{M}(\vect X)>m}$, where by varying the deterministic parameter $m$, the (complementary) distribution of $\mathcal{M}(\vect X)$ can be obtained. We rewrite Eq.\eqref{PerfectISLin}, the optimal importance density of an equivalent linear system, as
\begin{equation}\label{PerfectISProb}
h_L^*(\vect x)=\dfrac{\1{\mathcal{M}_L(\vect x)>m} f_{\vect X}(\vect x)}{\int_{\vect x\in\rn}\1{\mathcal{M}_L(\vect x)>m}f_{\vect X}(\vect x)\,d\vect x}=\dfrac{\1{\mathcal{M}_L(\vect x)>m}f_{\vect X}(\vect x)}{P_L}\,,
\end{equation}
where the normalizing constant $P_L$ is the probability $\Prob{\mathcal{M}_L(\vect X)>m}$. There is no guarantee that $h_L^*(\vect x)$ would dominate $\1{\mathcal{M}(\vect x)>m}f_{\vect X}(\vect x)$\footnote{Typically, the support $\lbrace\vect x:\1{\mathcal{M}_L(\vect x)>m}f_{\vect X}(\vect x)>0\rbrace$ of $h_L^*(\vect x)$ intersects $\lbrace\vect x:\1{\mathcal{M}(\vect x)>m}f_{\vect X}(\vect x)>0\rbrace$, i.e., there could be regions where $h_L^*(\vect x)$ is zero while $\1{\mathcal{M}(\vect x)>m}f_{\vect X}(\vect x)$ is not.}; therefore, $h_L^*(\vect x)$ cannot be used directly as an importance density to estimate $\E{Q}$. A straightforward remedy is to relax the hard indicator $\mathbb{I}$ into a smooth indicator $\hat{\mathbb{I}}$,
\begin{equation}\label{Eq:SmoothInd}
\hat{\mathbb{I}}(\mathcal{M}_L(\vect x)>m;\lambda)=\dfrac{1}{1+e^{\lambda(m-\mathcal{M}_L(\vect x))}}\,,
\end{equation}
where $\lambda>0$ is a parameter. For large $\lambda$, the smooth indicator closely resembles the hard indicator. Using Eq.\eqref{Eq:SmoothInd}, Eq.\eqref{PerfectISProb} is rewritten as
\begin{equation}\label{PerfectISProbS}
h_L^*(\vect x;\lambda)=\dfrac{\hat{\mathbb{I}}(\mathcal{M}_L(\vect x)>m;\lambda) f_{\vect X}(\vect x)}{\int_{\vect x\in\rn}\hat{\mathbb{I}}(\mathcal{M}_L(\vect x)>m;\lambda)f_{\vect X}(\vect x)\,d\vect x}=\dfrac{\hat{\mathbb{I}}(\mathcal{M}_L(\vect x)>m;\lambda)f_{\vect X}(\vect x)}{\hat P_L(\lambda)}\,.
\end{equation}
Since $\hat{\mathbb{I}}$ is nonzero everywhere, $h_L^*(\vect x;\lambda)$ is a valid importance density to estimate $\E{Q}$. The parameter $\lambda$ can be tuned to be as small as possible, under the constraint that if $h_L^*(\vect x;\lambda)$ is used to estimate $P_L$ of the linear system, the efficiency should reach a specified threshold. This tuning aims to address the trade-off between sampling bias and sampling variance of the importance density. For a large $\lambda$, the importance density is concentrated, but the estimation will be biased.

Eq.\eqref{PerfectISProbS} provides a feasible adaptation of the original framework to estimate rare event probabilities. Alternatively, it is possible to use Eq.\eqref{PerfectISProb} in its original form and still develop a theoretically correct approach to estimate $\E{Q}$. Specifically, we construct the following equation/identity
\begin{equation}\label{CondProb}
\begin{aligned}
\E{Q}&=P_L\dfrac{\int_{\vect x\in\rn}\1{\mathcal{M}(\vect x)>m}\1{\mathcal{M}_L(\vect x)>m}h^*_L(\vect x)\,d\vect x}{\int_{\vect x\in\rn}\1{\mathcal{M}(\vect x)>m}\1{\mathcal{M}_L(\vect x)>m}h^*(\vect x)\,d\vect x}\\
&\approx P_L\dfrac{\left\langle\1{\mathcal{M}(\vect x)>m}\1{\mathcal{M}_L(\vect x)>m}\right\rangle_{h^*_L}}{\left\langle\1{\mathcal{M}(\vect x)>m}\1{\mathcal{M}_L(\vect x)>m}\right\rangle_{h^*}}\\
&=P_L\dfrac{\left\langle\1{\mathcal{M}(\vect x)>m}\right\rangle_{h^*_L}}{\left\langle\1{\mathcal{M}_L(\vect x)>m}\right\rangle_{h^*}}\,.
\end{aligned}
\end{equation}
It is straightforward to verify Eq.\eqref{CondProb} using the definitions of $h^*_L(\vect x)$ and $h^*(\vect x)$. The numerator and denominator can be estimated by sampling from $h^*_L(\vect x)$ and $h^*(\vect x)$, respectively. If the linear system response somewhat resembles the nonlinear system response, the probabilities in the numerator and denominator would be relatively large. It is worth mentioning that similar to Eq.\eqref{ISLin}, Eq.\eqref{CondProb} can also be interpreted as applying a correction factor to the equivalent linear system solution. Moreover, Eq.\eqref{CondProb} can be viewed as a variant of the \textit{Bennett acceptance ratio
}, which was originally introduced in \cite{bennett1976efficient} and recently studied in the context of rare event simulation in \cite{karl}. 

Finally, it is {worth} mentioning that the above two approaches of adapting $h_L^*(\vect x)$ for rare-event simulations are by no means unique. For example, in the multi-fidelity UQ framework \cite{peherstorfer2018survey}, it is suggested to fit a Gaussian mixture model over samples of $h_L^*(\vect x)$, and then use the Gaussian mixture as the importance density. Unfortunately, this approach does not work in our context because we aim to solve random vibration problems with high-dimensional uncertainties, and Gaussian mixture-based importance sampling does not scale well with dimensionality (see \cite{WANG201642} for an in-depth discussion). It might be promising to investigate the use of directional distribution models (which scale better with dimensionality than Gaussian) such as a mixture of von Mises-Fisher distributions. However, notice that nonlinear random vibration problems such as the first-passage probability estimation can have a fairly complex integration domain; fitting a simple parametric distribution model over a complex high-dimensional domain can be challenging. 

\subsection{Optimization}
\noindent For rare event probability estimations, the iterative approach developed in Section \ref{Sec:OOIS} can be too fragile to be effective. This is because with some bad luck, almost all samples generated from $h_L^*(\vect x;\vect\theta_j^*)$ of a trial linear system can lead to $\1{\mathcal{M}(\vect x)>m}$ being zero. In principle, this issue can be mitigated by performing relaxation on the indicator function $\1{\mathcal{M}(\vect x)>m}$, e.g., via a smooth indicator $\hat{\mathbb{I}}(\mathcal{M}(\vect x)>m; \lambda)$, or via introducing $\lambda>0$ into $\1{\mathcal{M}(\vect x)+\lambda>m}$. However, to make the optimization more robust, we consider a direct approach. Specifically, instead of iterating with trial linear systems, we perform a one-step optimization using MCMC samples from $h^*(\vect x)\propto\1{\mathcal{M}(\vect x)>m}f_{\vect X}(\vect x)$. Moreover, instead of attempting to reproduce the topological quantity $\1{\mathcal{M}(\vect x)>m}$, which contains information far less than the model $\mathcal{M}(\vect x)$, we let the optimized linear system imitate $\mathcal{M}(\vect x)$ at the localized region $\lbrace\vect x:\1{\mathcal{M}(\vect x)>m}=1\rbrace$. Consequently, following a similar derivation of Section \ref{Sec:OOIS}, the variance minimization-based optimization is
\begin{equation}\label{OptProbMV}
\begin{aligned}
\vect\theta^*&=\mathop{\arg\min}_{\vect\theta}\left\lbrace\int_{\vect x\in\rn}\Big|\dfrac{\mathcal{M}(\vect x)}{\mathcal{M}_L(\vect x;\vect\theta)}\Big|h^*(\vect x)\,d\vect x\,\int_{\vect x\in\rn}\Big|\dfrac{\mathcal{M}_L(\vect x;\vect\theta)}{\mathcal{M}(\vect x)}\Big|h^*(\vect x)\,d\vect x\right\rbrace\\
&\approx\mathop{\arg\min}_{\vect\theta}\left\lbrace\left\langle\Big|\dfrac{\mathcal{M}(\vect x)}{\mathcal{M}_L(\vect x;\vect\theta)}\Big|\right\rangle_{h^*}\left\langle\Big|\dfrac{\mathcal{M}_L(\vect x;\vect\theta)}{\mathcal{M}(\vect x)}\Big|\right\rangle_{h^*}\right\rbrace\,,
\end{aligned}
\end{equation}
The cross entropy-based optimization is
\begin{equation}\label{OptProbCE}
\begin{aligned}
\vect\theta^*&=\mathop{\arg\max}_{\vect\theta}\int_{\vect x\in\rn}\log\dfrac{|\mathcal{M}_L(\vect x;\vect\theta)|}{\mu_{|L|}(\vect\theta)}h^*(\vect x)\,d\vect x\\
&\approx\mathop{\arg\max}_{\vect\theta}\left\langle\log|\mathcal{M}_L(\vect x;\vect\theta)|-\log\Big\langle\Big|\dfrac{\mathcal{M}_L(\vect x;\vect\theta)}{\mathcal{M}(\vect x)}\Big|\Big\rangle_{h^*}\right\rangle_{h^*}\,.
\end{aligned}
\end{equation}
Finally, it is important to restate that Eq.\eqref{OptProbMV} and Eq.\eqref{OptProbCE} are typically invariant with respect to the scaling factor $b$ of the linear system; because these equations try to establish a linear system with the response being maximally correlated to the nonlinear system response. However, the scaling factor $b$ is non-negligible when the indicator function $\1{\mathcal{M}_L(\vect x;\vect\theta)>m}$ is of concern. Therefore, after solving Eq.\eqref{OptProbMV} or Eq.\eqref{OptProbCE} for the linear system parameters other than $b$, $b$ can then be fixed to
\begin{equation}
b^*=\mathop{\arg\min}_{b}\left(\left\langle\mathcal{M}_L(\vect x;b)\right\rangle_{h^*}-\left\langle\mathcal{M}(\vect x)\right\rangle_{h^*}\right)^2\,.
\end{equation}
This optimization does not require additional calls of the nonlinear model because the samples of $h^*(\vect x)$ can be inherited from the previous step of solving Eq.\eqref{OptProbMV} or Eq.\eqref{OptProbCE}. The implementation details for the approach based on indicator relaxation (Eq.\eqref{PerfectISProbS}) and the approach based on conditioning (Eq.\eqref{CondProb}) are { developed in \ref{Fig:ROIS} and \ref{Fig:COIS}}, respectively.

\section{The General Framework}
\begin{figure}[H]
    \centering
    \includegraphics[scale=0.6]{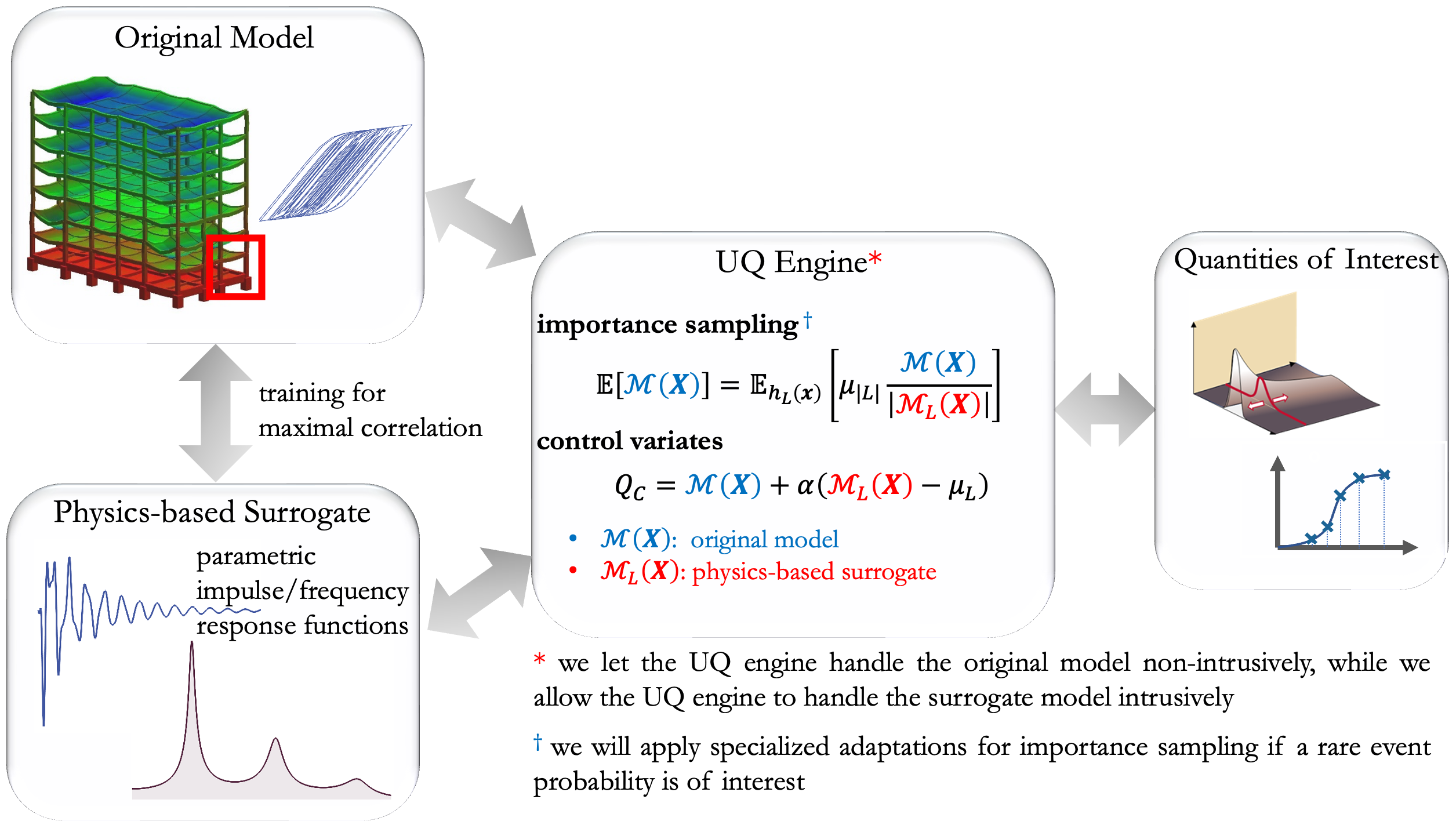}
    \caption{{\textbf{General framework of the optimized equivalent linearization method}}. \textit{The proposed optimized equivalent linearization method can be reworded as a linear dynamic system-based surrogate modeling approach for nonlinear stochastic dynamics analysis. The linear system is represented by a parametric impulse/frequency response function with physical parameters such as damping, stiffness, modal participation factor, degree-of-freedom, etc. This is a parsimonious surrogate model for a high-dimensional uncertainty quantification problem. The computational flows (arrows) are bidirectional since the modules are coupled, e.g., the quantity of interest evaluations can guide computational efforts in the surrogate model training and UQ solver}.}
    \label{Fig:UQ}
\end{figure}

\noindent Figure \ref{Fig:UQ} schematically summarizes the proposed optimized equivalent linearization method. 
It is apparent from the figure that the linear system model can be replaced by some parametric nonlinear system model, and the framework is not restricted to nonlinear random vibration problems. The critical feature of the framework is that the adopted surrogate model should be physics-based and thus interpretable to engineers; therefore, the vision of the classical equivalent linearization methodology is preserved.

\section{Numerical Examples}
\subsection{Stochastic Excitation Model}
\noindent For all subsequent examples, we consider a stochastic acceleration process $a(t)$ as the external excitation. The acceleration process is described by the following power spectrum density model suggested by Clough and Penzien \cite{clough1975st}, which is widely used in earthquake engineering.  
\begin{equation}\label{KTSpec}
S_g(\omega)=S_0\frac{\omega_f^4+4\zeta_f^2\omega_f^2\omega^2}{(\omega_f^2-\omega^2)^2+4\zeta_f^2\omega_f^2\omega^2}\dfrac{\omega^4}{(\omega_s^2-\omega^2)^2+4\zeta_s^2\omega_s^2\omega^2}\,,
\end{equation}
where $S_0$ is a scale factor, $\omega_f=15\,\rm{rad/s}$ and $\zeta_f=0.6$ are the filter parameters representing, respectively, the natural frequency and damping ratio of the soil layer, and $\omega_s=1.5\,\rm{rad/s}$ and $\zeta_s=0.6$ are parameters of a second filter introduced to assure finite variance of the displacement. The duration of the process is set to 10 seconds. The parameters are selected such that the process is wide-band, i.e., the randomness cannot be significantly reduced to make the random vibration problem easier. The random realization of $a(t)$ is represented by \cite{shinozuka1972digital}\cite{chatfield2013analysis}
\begin{equation}
a(t)=\sum_{i=1}^{n/2}\sigma_i\left(x_{2i-1}\sin\omega_it+x_{2i}\cos\omega_it\right)\,,
\end{equation}
where the frequency domain $[0,\omega_{n/2}]$ is uniformly discretized  with the increment $\Delta\omega$, $\sigma_i=\sqrt{2S_g(\omega_i)\Delta\omega}$, and $x_i$ are realizations of independent standard Gaussian random variables. We set $n=200$ and $\omega_{n/2}=15\pi\,\rm{rad/s}$. Therefore, eventually, the stochastic process is described in the probability space of $\vect X=[X_1,X_2,\dots,X_{200}]$ with 200 Gaussian random variables. 

\subsection{Optimization and Sampling}
\noindent The optimization problems of identifying the equivalent linear system parameters are solved by the interior-point method \cite{byrd2000trust}\cite{waltz2006interior} with a multistart setting \citep{ugray2007scatter}. Specifically, the number of random starting points for the multistart optimization is set to be proportional to $DoF$, i.e., $\alpha\cdot DoF$, where $\alpha=3$ seems sufficient to obtain stable solutions for the studied examples. For the adaptive importance sampling-enhanced ELM, there is an additional variability in choosing the optimization objective: variance minimization or cross-entropy minimization. However, we did not observe any convincing evidence to favor one over the other; therefore, one could implement either.

Except for probability estimations, the Gibbs sampling \cite{andrieu2003introduction} with a uniform proposal distribution is adopted to sample from the equivalent linear system-based importance density $h_L^*(\vect x)$. The step size of the Gibbs sampler is adaptively tuned using acceptance rates of Markov chains \cite{haario2001adaptive}\cite{WANG2019}. For (rare event) probability estimations, the rejection sampling-based Hamiltonian Monte Carlo method \cite{WANG2019} is adopted to sample from the nonlinear system and linear system-based importance densities. To reduce the number of nonlinear system simulations, we let the initial samples of the Hamiltonian Monte Carlo sampler be located in the critical region that defines the rare event. This is done via a secant method to project random samples of $f_{\vect X}(\vect x)$ to the critical region. %For readers familiar with MCMC techniques, note that we did not use Hamiltonian Monte Carlo to sample from $h_L^*(\vect x)$ for problems other than probability estimations. This is because Hamiltonian Monte Carlo requires gradients of the distribution function, and $h_L^*(\vect x)$ can be non-differentiable (with respect to $\vect x$) for some QoI definitions (e.g., those involving $\max$ or $\min$ operators). On the other hand, in the context of probability estimations, we can use rejection sampling-based Hamiltonian Monte Carlo because rejection sampling effectively eliminates the requirement of $h^*(\vect x)$ or $h_L^*(\vect x)$ being differentiable; it only requires $f_{\vect X}(\vect x)$ being differentiable. 

\subsection{Performance Evaluations}
\noindent The performance of the proposed optimized equivalent linearization methodology is investigated by comparing with solutions of the direct Monte Carlo simulation. In the following numerical examples, we push the random vibration problems to strong nonlinearity regions such that the peak response can be 12 times larger than the yield point; we also estimate rare event probabilities (first-passage probabilities) as low as $10^{-8}$. In such a problem setup, the solutions of conventional ELMs are not comparable. Moreover, popular surrogate modeling techniques such as the active learning-based Kriging method \cite{echard2011ak} are not considered for comparison because they cannot solve nonlinear random vibration problems with high-dimensional uncertainties. 

\subsection{Cubic Oscillator}
\begin{figure}[H]
    \centering
    \includegraphics[scale=0.47]{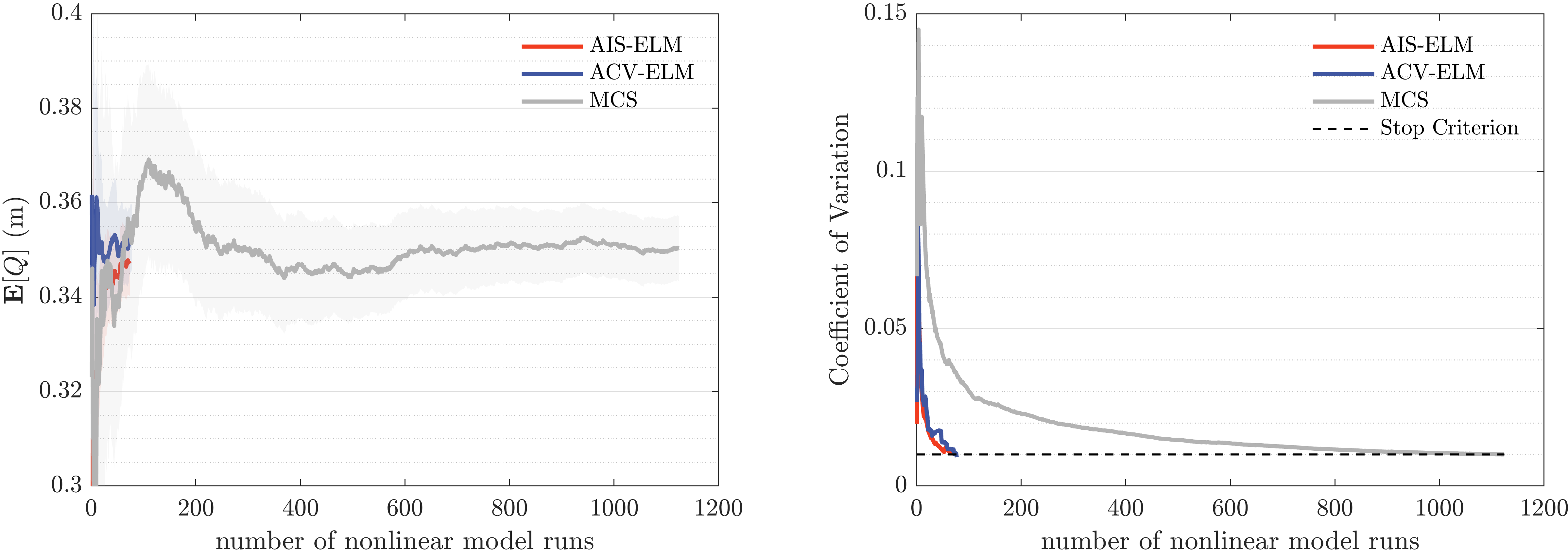}
    \caption{\textbf{Simulation result of the cubic oscillator for the mean peak absolute response}. \textit{To achieve the same coefficient of variation of $1\%$, the ACV-ELM requires $50$(for optimization)+$78=128$ nonlinear model runs, and the estimate is $0.35\,\rm{m}$; the AIS-ELM requires $50$(for optimization)$+74=124$ nonlinear model runs, and the estimate is $0.35\,\rm{m}$; the direct Monte Carlo simulation requires $1124$ nonlinear model runs, and the estimate is $0.35\,\rm{m}$. The shaded areas represent approximate confidence intervals obtained from the mean plus and minus $1.96$ standard deviations}.}
    \label{Fig:RunCubic}
\end{figure}

\begin{figure}[H]
    \centering
    \includegraphics[scale=0.47]{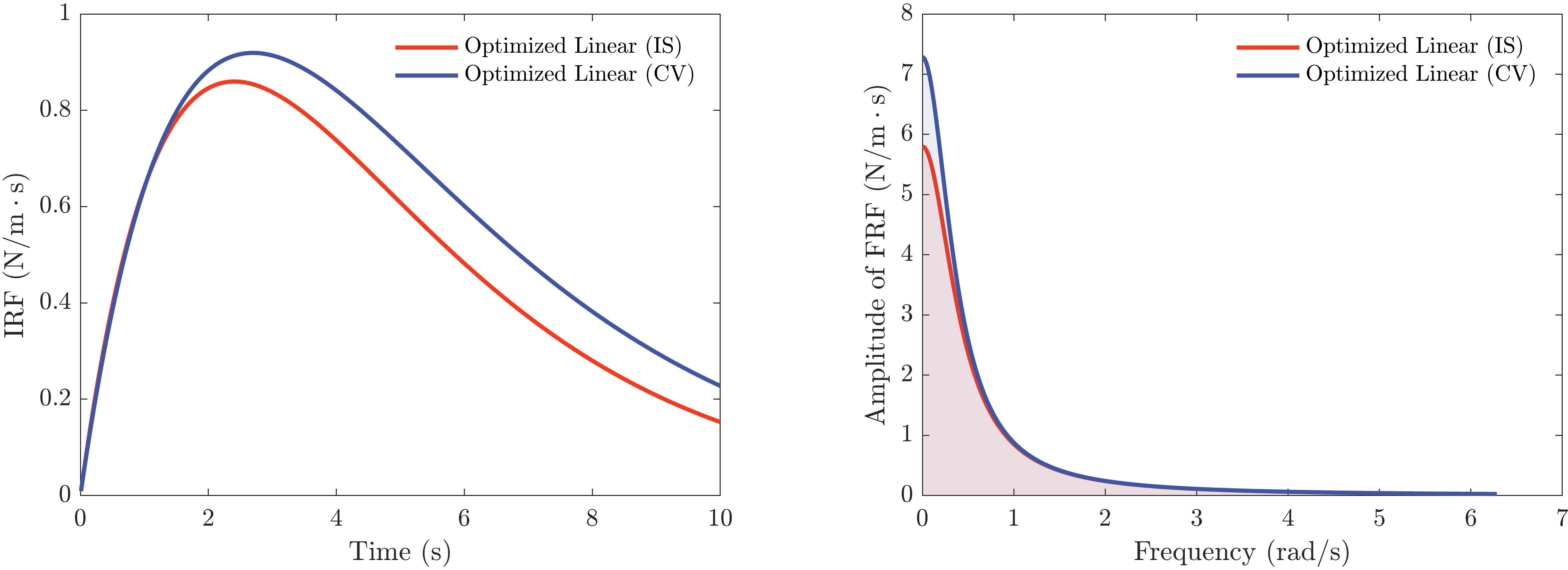}
    \caption{\textbf{Impulse and frequency response functions of the optimized linear system for the mean peak response of the cubic oscillator}. \textit{The impulse/frequency response functions obtained from both methods exhibit dominant low-frequency components}.}
    \label{Fig:FRFCubic}
\end{figure}

\noindent Consider a single degree-of-freedom cubic oscillator described by the equation of motion
\begin{equation}
 {\ddot z}(t)+{\dot z}(t)+z^3(t)=-a(t)\,.
\end{equation}
Setting $S_0=0.03\,\rm{m^2/s^3}$\footnote{We solved this example multiple times by picking $S_0$ values from the range $[0,1]$. We concluded that a high correlation between optimized linear and nonlinear system responses is achievable regardless of the specific $S_0$ value (within $[0,1]$). Notice that the range $[0,1]$ is wide enough to generate strong nonlinearities.}, the QoI is defined as the peak absolute deformation so that the $\E{Q}$ to be estimated is
\begin{equation}\label{peak}
\E{Q}=\E{\sup_{t\in[0,10]}\left\vert Z(t)\right\vert}\,.
\end{equation}
The Adaptive Control Variate-enhanced ELM (ACV-ELM) and the Adaptive Importance Sampling-enhanced ELM (AIS-ELM) are used to compute $\E{Q}$. The target coefficient of variation is set to $1\%$ for all methods considered. A simulation result is reported in Figure \ref{Fig:RunCubic}. The impulse response and frequency response functions of the optimized linear systems are illustrated in Figure \ref{Fig:FRFCubic}.

\begin{figure}[H]
    \centering
    \includegraphics[scale=0.45]{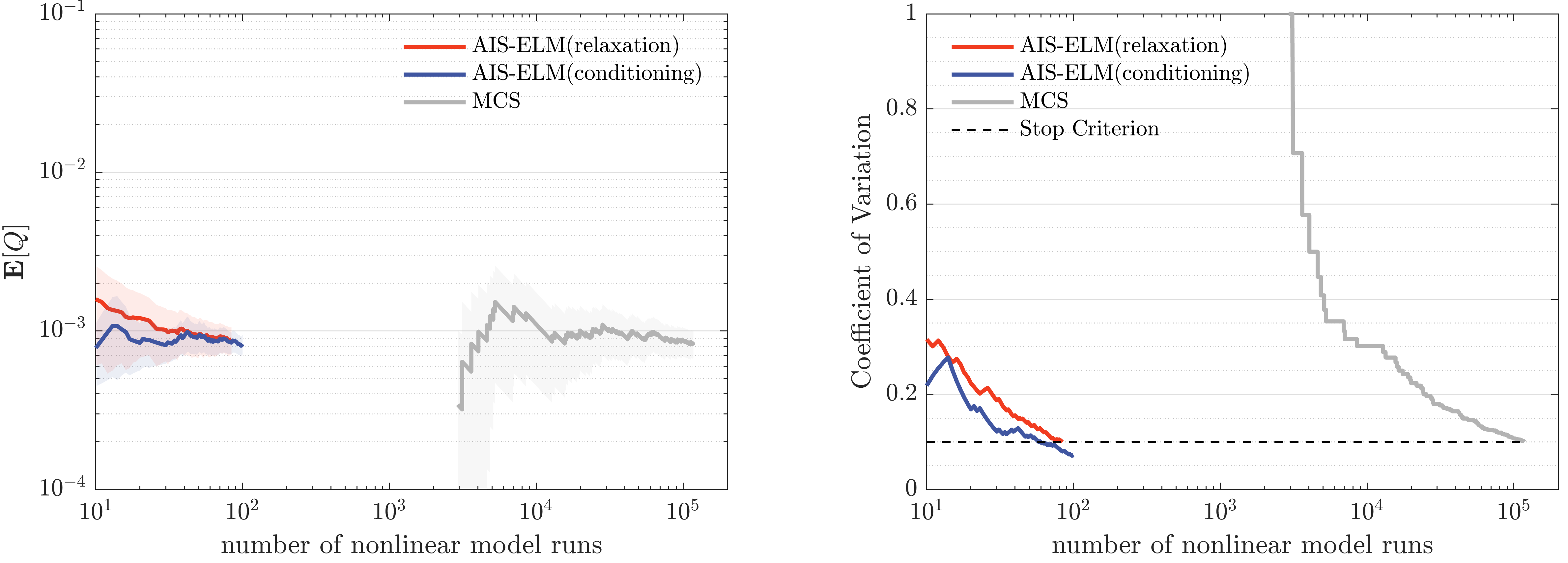}
    \caption{\textbf{Simulation result for the first-passage probability estimation of the cubic oscillator}. \textit{To achieve the same coefficient of variation of $10\%$, the AIS-ELM based on relaxation requires $653$(for optimization)$+84=737$ nonlinear model runs, and the probability estimate is $8.84\times10^{-4}$; the AIS-ELM based on conditioning requires $653$(for optimization)$+500=1153$ nonlinear model runs, and the probability estimate is $8.28\times10^{-4}$; the direct Monte Carlo simulation requires $1.18\times10^5$ nonlinear model runs, and the probability estimate is $8.50\times10^{-4}$. The AIS-ELM based on conditioning involves ``hidden" nonlinear model runs embedded in the MCMC sampler, which cannot be shown in the figure (but are already included in the reported computational cost). The shaded areas represent approximate confidence intervals obtained from the mean plus and minus $1.96$ standard deviations}.}
    \label{Fig:RunFPCubic}
\end{figure}

Next, we consider the first-passage probability of the peak absolute displacement, i.e.,
\begin{equation}
\E{Q}=\E{\1{\sup_{t\in[0,10]}\left\vert Z(t)\right\vert>m}}=\Prob{\sup_{t\in[0,10]}\left\vert Z(t)\right\vert>m}\,,
\end{equation}
where we set the response threshold $m$ to $0.8\,\rm{m}$. The AIS-ELM approaches based on relaxation and conditioning are used to compute the first-passage probability. The target coefficient of variation is set to $10\%$ for all methods considered. A simulation result is reported in Figure \ref{Fig:RunFPCubic}.

\begin{figure}[H]
    \centering
    \includegraphics[scale=0.41]{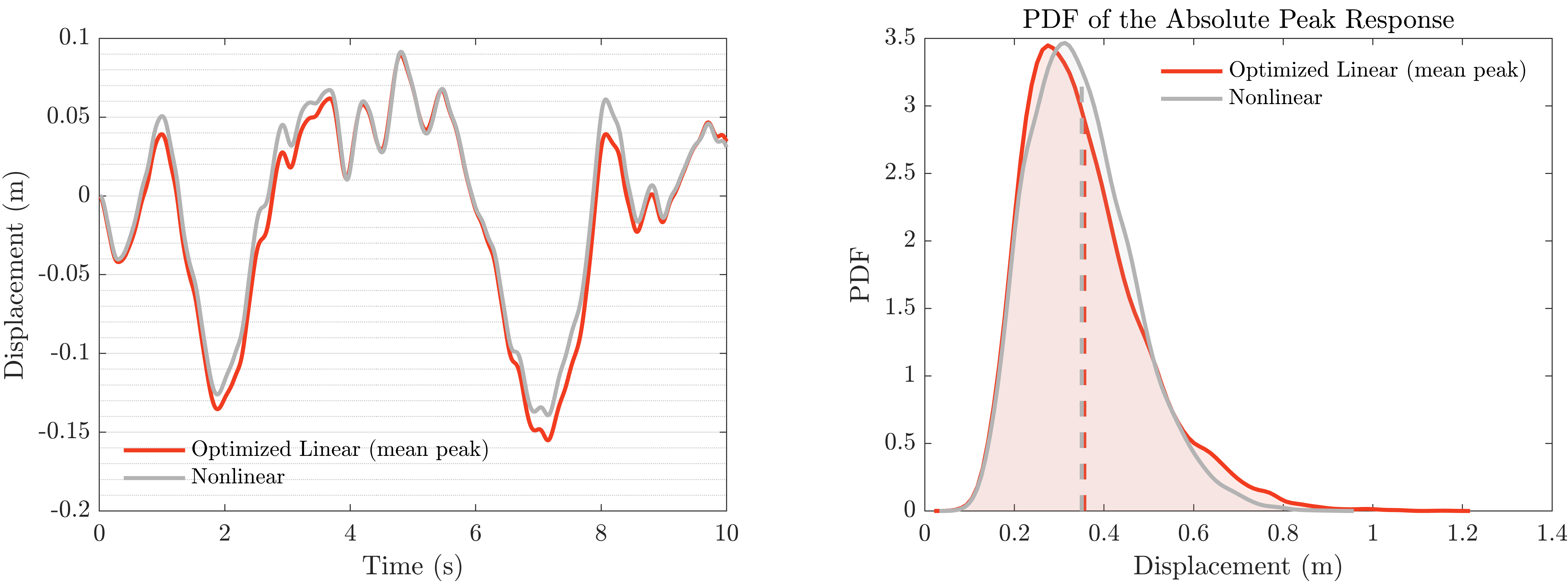}
    \caption{\textbf{Response trajectories and distributions of the optimized linear system}. \textit{The  linear system could imitate the actual trajectory of the cubic oscillator (although the optimization objective is to match the mean peak response); their peak response distributions slightly differ. The mean peak responses (dashed lines) of the two systems match, because the linear system is optimized to estimate the mean peak response}.}
    \label{Fig:TSCubic}
\end{figure}

The results from mean peak response and first-passage probability estimations clearly demonstrate the efficiency and accuracy of the optimized ELM. Without compromising the accuracy, the proposed method can be two orders more efficient than the direct MCS for first-passage probability estimations. Finally, the left panel of Figure \ref{Fig:TSCubic} illustrates the response trajectories of the linear systems (optimized for mean peak response estimation) compared with the actual response of the cubic oscillator under one random realization of the ground motion; the right panel compares the peak response distributions.

\subsection{Bouc-wen Oscillator}
\noindent Consider a single degree-of-freedom oscillator with Bouc-Wen hysteretic constitutive model. The equation of motion is expressed as
\begin{equation}
\begin{array}{lr}\label{BW}
  {\ddot z}(t)+2\omega_n\zeta_n{\dot z}(t)+\omega_n^2(\alpha_{bw} z(t)+(1-\alpha_{bw})z_{bw}(t))=-a(t) \\
  \dot{z}_{bw}(t)=-\gamma\left|\dot{z}(t)\right|\left|z_{bw}(t)\right|^{\bar{n}-1}z_{bw}(t)-\eta\left|z_{bw}(t)\right|^{\bar n}\dot{z}(t)+A\dot{z}(t)\,,
\end{array}
\end{equation}
where we set $\omega_n=10\sqrt{10}\,\rm{rad/s}$, $\zeta_n=0.05$, $\alpha_{bw}=0.1$, $\bar{n}=5$, $A=1$, $\gamma=\eta=1/(2z_y^{\bar{n}})$, where $z_y=0.01\,\rm{m}$ is the effective yield displacement. 

Setting $S_0=0.03\,\rm{m^2/s^3}$, the proposed method is applied to estimate the mean peak response and first-passage probability for a threshold of $0.08\,\rm{m}$. For all methods considered, the target coefficient of variation is $1\%$ for mean peak response and $10\%$ for first-passage probability. The results are reported in Table \ref{tab1}.

\begin{table}[H]
  \caption{\textbf{Mean peak response and first-passage probability estimations}.}
  \label{tab1}
  \centering
  
  \begin{tabular}{c c c c c}
    \toprule
    Methods & \multicolumn{2}{c}{Mean absolute peak} & \multicolumn{2}{c}{First-passage probability} \\
   & $\E{Q}$&$N$& $\E{Q}$&$N$\\
    \midrule
   ACV-ELM & $3.18\times10^{-2}$&219&-&-\\
   AIS-ELM & $3.15\times10^{-2}$ &349&-&-\\
   AIS-ELM (conditioning) &-&- &$1.36\times10^{-4}$&$1.30\times10^{3}$\\
   AIS-ELM (relaxation) &-&-  &$1.28\times10^{-4}$&$1.30\times10^{3}$\\
   MCS & $3.18\times10^{-2}$&584&$1.37\times10^{-4}$&$7.28\times10^{5}$ \\
    \bottomrule
  \end{tabular}
\end{table}

\begin{figure}[H]
    \centering
    \includegraphics[scale=0.5]{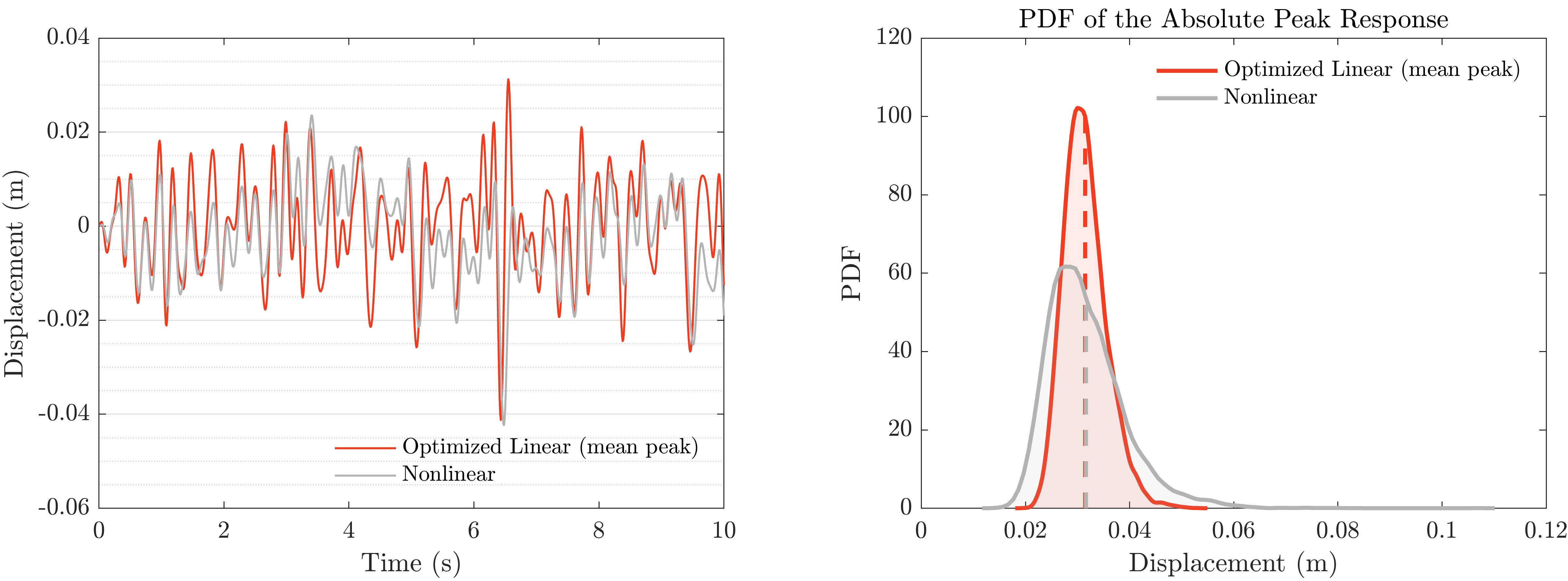}
    \caption{\textbf{Response trajectories and distributions of the optimized linear system}. \textit{The  linear system response trajectory mildly resembles the nonlinear one. The peak response distributions noticeably differ. The mean peak responses (dashed lines) of the two systems match, because the linear system is optimized to estimate the mean peak response}.}
    \label{Fig:TSBoucWen}
\end{figure}

Similar to the previous example, the proposed method can be two orders more efficient than the direct MCS for first-passage probability estimations. Figure \ref{Fig:TSBoucWen} illustrates the response trajectories of the linear systems (optimized for mean peak response estimation) compared with the nonlinear response, under one random realization of the ground motion; the right panel compares the peak response distributions. 

Finally, using this example, we demonstrate that the optimized ELM readily applies to nonstationary problems without requiring any modification. This is because, regardless of the stochastic excitation model, the stochastic process samples are propagated into pre-specified response quantities to optimize parametric linear systems. A piecewise temporal modulating function 
\begin{equation}
m(t)=\left\lbrace\begin{aligned}
&(t/T_1)^2\,,&&t\in[0,T_1]\\
&1\,,&&t\in(T_1,T_2)\\
&e^{-(t-T_2)^2}\,,&&t\geq T_2
\end{aligned}\right.
\end{equation}
where $T_1=4$ and $T_2=7$, is applied to the stationary excitation model with $S_0=0.2\,\rm{m^2/s^3}$. The quantities in Table \ref{tab1} are recomputed to generate Table \ref{tab2}. The results suggest significant efficiency gains for using the optimized ELM in nonstationary nonlinear random vibration analysis. 

\begin{table}[H]
  \caption{\textbf{Mean peak response and first-passage probability estimates for nonstationary excitation}.}
  \label{tab2}
  \centering
  
  \begin{tabular}{c c c c c}
    \toprule
    Methods & \multicolumn{2}{c}{Mean absolute peak} & \multicolumn{2}{c}{First-passage probability} \\
   & $\E{Q}$&$N$& $\E{Q}$&$N$\\
    \midrule
   ACV-ELM & $3.22\times10^{-2}$&367&-&-\\
   AIS-ELM & $3.19\times10^{-2}$ &440&-&-\\
   AIS-ELM (conditioning) &-&- &$1.10\times10^{-3}$&$1.30\times10^{3}$\\
   AIS-ELM (relaxation) &-&-  &$1.06\times10^{-3}$&$0.85\times10^{3}$\\
   MCS & $3.24\times10^{-2}$&844&$0.97\times10^{-3}$&$1.03\times10^{5}$ \\
    \bottomrule
  \end{tabular}
\end{table}

\subsection{Multiple Degree of Freedom System}
\begin{figure}[H]
    \centering
    \includegraphics[scale=0.75]{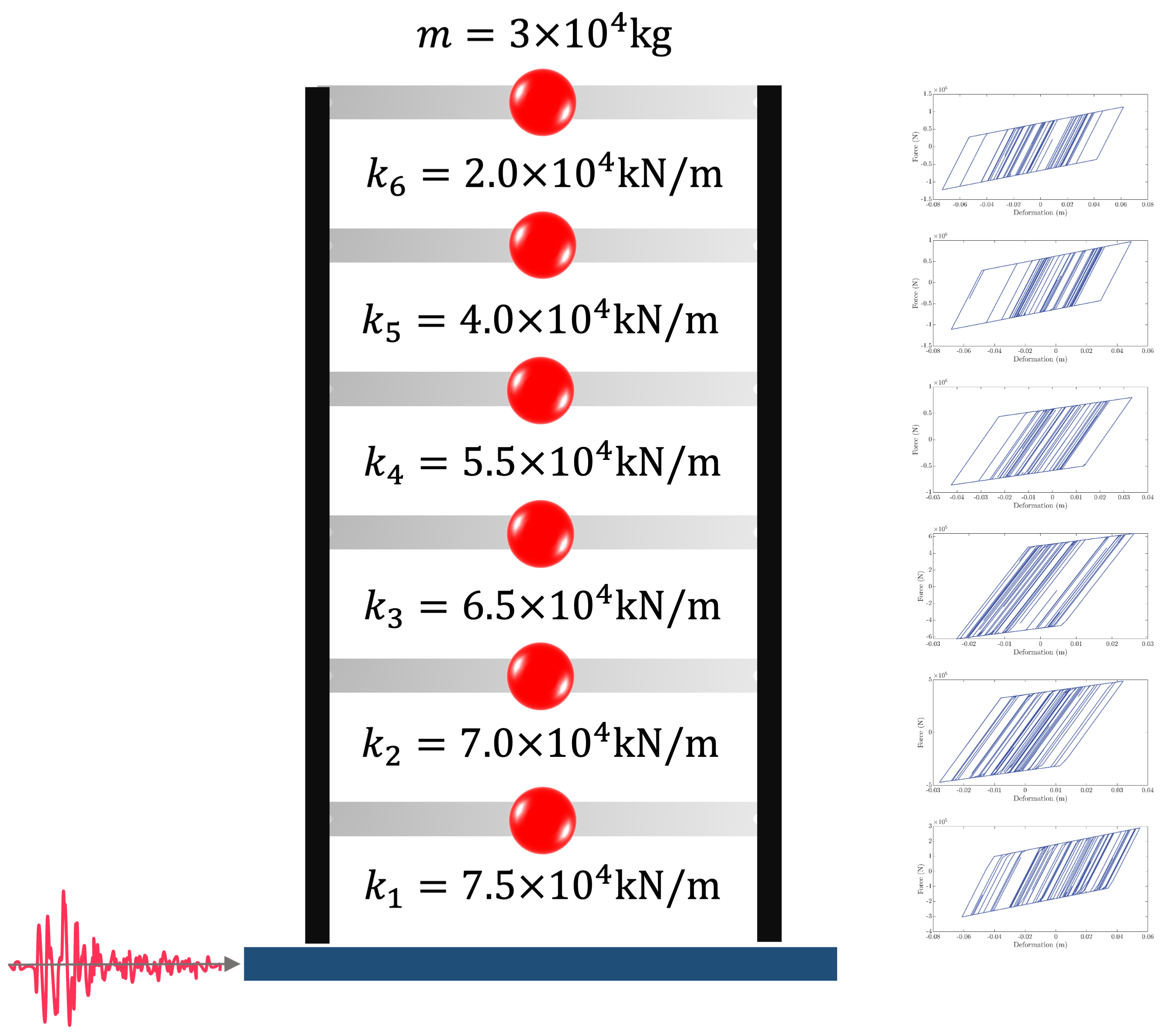}
    \caption{\textbf{Shear-building model and the force-deformation behavior for each story}.}
    \label{Fig:ShearWall}
\end{figure}

\noindent Consider a 6-Degree-of-Freedom shear-building model \cite{fujimura2007tail}\cite{wang2017equivalent} shown in Figure \ref{Fig:ShearWall}. The constitutive relationship between force and inter-story deformation of each story is described by a bilinear hysteretic model shown in the figure. The yield deformation of each story is set to $0.01\,\rm{m}$. The structure has an initial fundamental period of $0.58$ seconds and a second mode period of $0.24$ seconds. Modal damping with a $5\%$ damping ratio for each mode is assumed. 

\begin{figure}[H]
    \centering
    \includegraphics[scale=0.6]{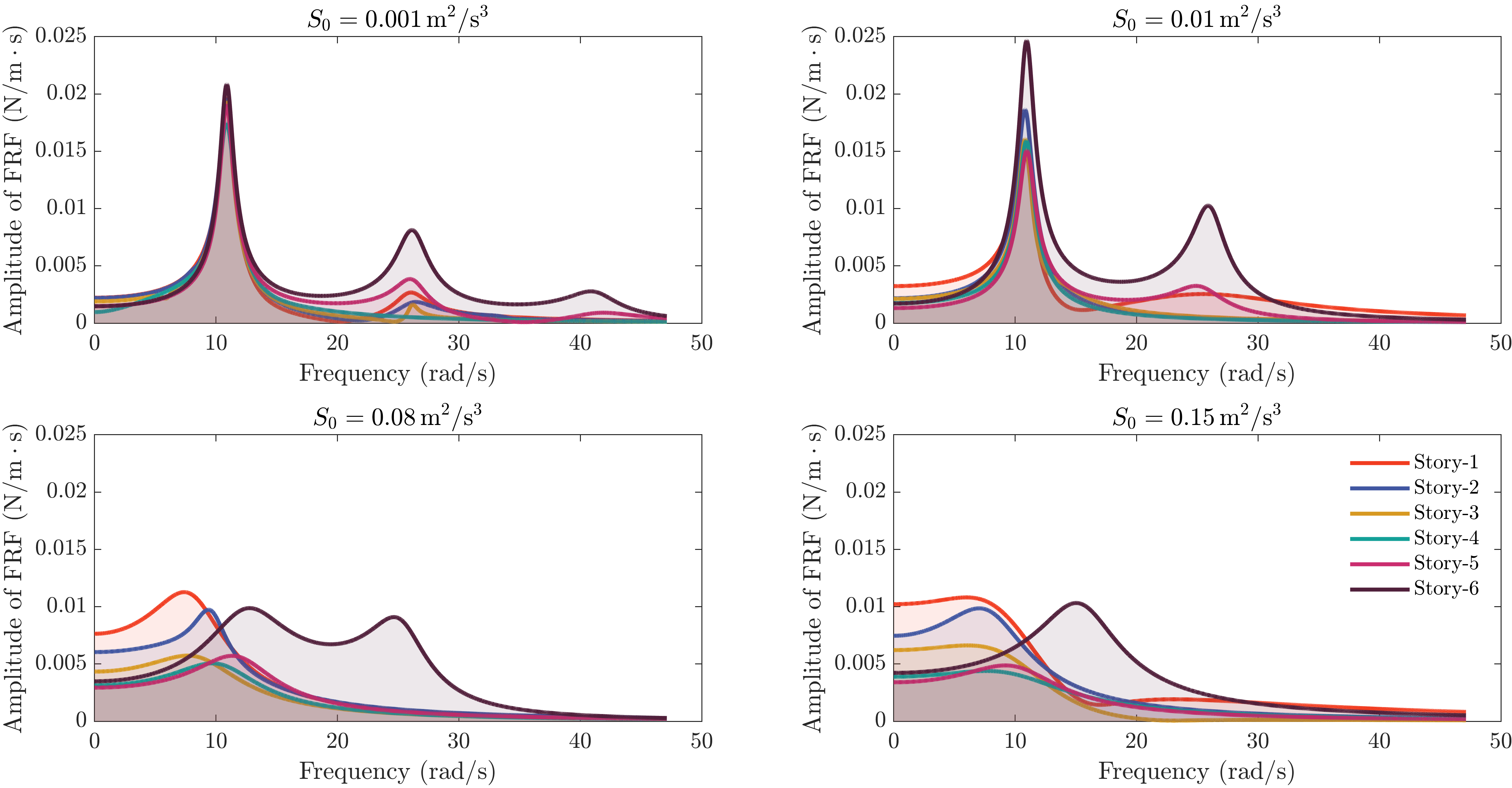}
    \caption{\textbf{Frequency response functions of the optimized linear systems for the mean peak response of different stories and excitation intensities}. \textit{The optimized linear system evolves with the degree of nonlinearity. At $S_0=0.001\,\rm{m^2/s^3}$, where the nonlinearity in the original system is not triggered, the frequency response functions of the optimized linear system are the same as that of the tangential linear system (at zero deformation), and the frequency response functions exhibit multimodal behavior. As the excitation intensity grows, the higher-frequency modes decay, and the lower-frequency component becomes dominant. The optimized $DoF$ values for equivalent linear systems are consistent with the above trend. For small $S_0$, the $DoF$ can be $4$ or $5$. As $S_0$ increases, the $DoF$ reduces into $1$ or $2$}.}
    \label{Fig:FRFDoF}
\end{figure}

It is expected that for different QoIs and different excitation intensities, the optimized linear systems differ. Figure \ref{Fig:FRFDoF} illustrates the frequency response functions of the optimized linear systems for estimating the mean peak absolute deformation of different stories under increasing excitation intensities. The result is obtained from the ACV-ELM, while  AIS-ELM would yield similar trend. 

\begin{figure}[H]
    \centering
    \includegraphics[scale=0.5]{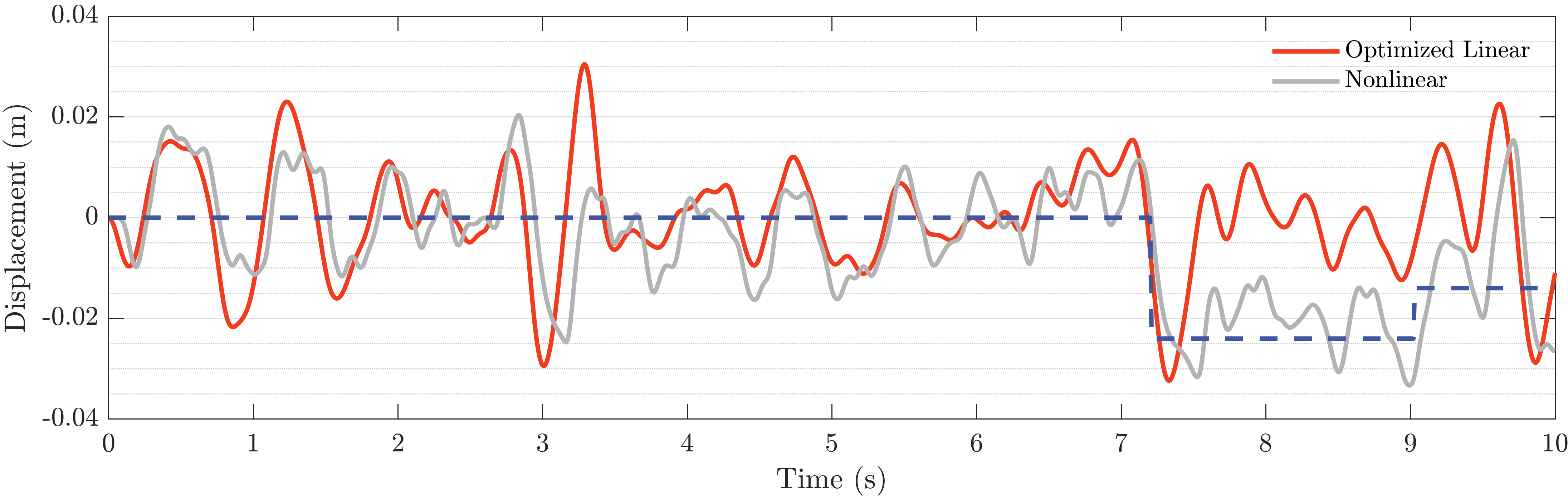}
    \caption{\textbf{Responses at the 5th story from the optimized linear and nonlinear systems}. \textit{There is a step-function-like trend (qualitatively illustrated as a dashed line) in the nonlinear system response, which cannot be modeled by a superposition of linear systems. The trend might be caused by the phase difference in plasticity development between consecutive floors. Notice that the dashed line is for qualitative illustration purpose only; trend detection algorithms for nonlinear response trajectories are out of our current scope}.}
    \label{Fig:TSMDoF}
\end{figure}

At $S_0=0.001\,\rm{m^2/s^3}$, where the nonlinearity in the original system is not triggered, the correlation coefficients between the peak responses of the optimized linear and original systems are around $99.99\%$. This result is theoretically trivial/obvious, but practically it implies that our algorithmic implementations are robust. At $S_0=0.15\,\rm{m^2/s^3}$, where the inter-story deformations can be $5$ times larger than the yield displacement, the correlation coefficients between the peak responses of the optimized linear and nonlinear systems are around $80\%$ for the 1st to 3rd story and can be as low as $40\%$ for the 5th story. We believe this low correlation coefficient cannot be improved without changing the design of the parametric linear system. To further investigate this issue, Figure \ref{Fig:TSMDoF} compares typical responses at the 5th story obtained from the optimized linear and nonlinear systems. 

\begin{figure}[H]
    \centering
    \includegraphics[scale=0.6]{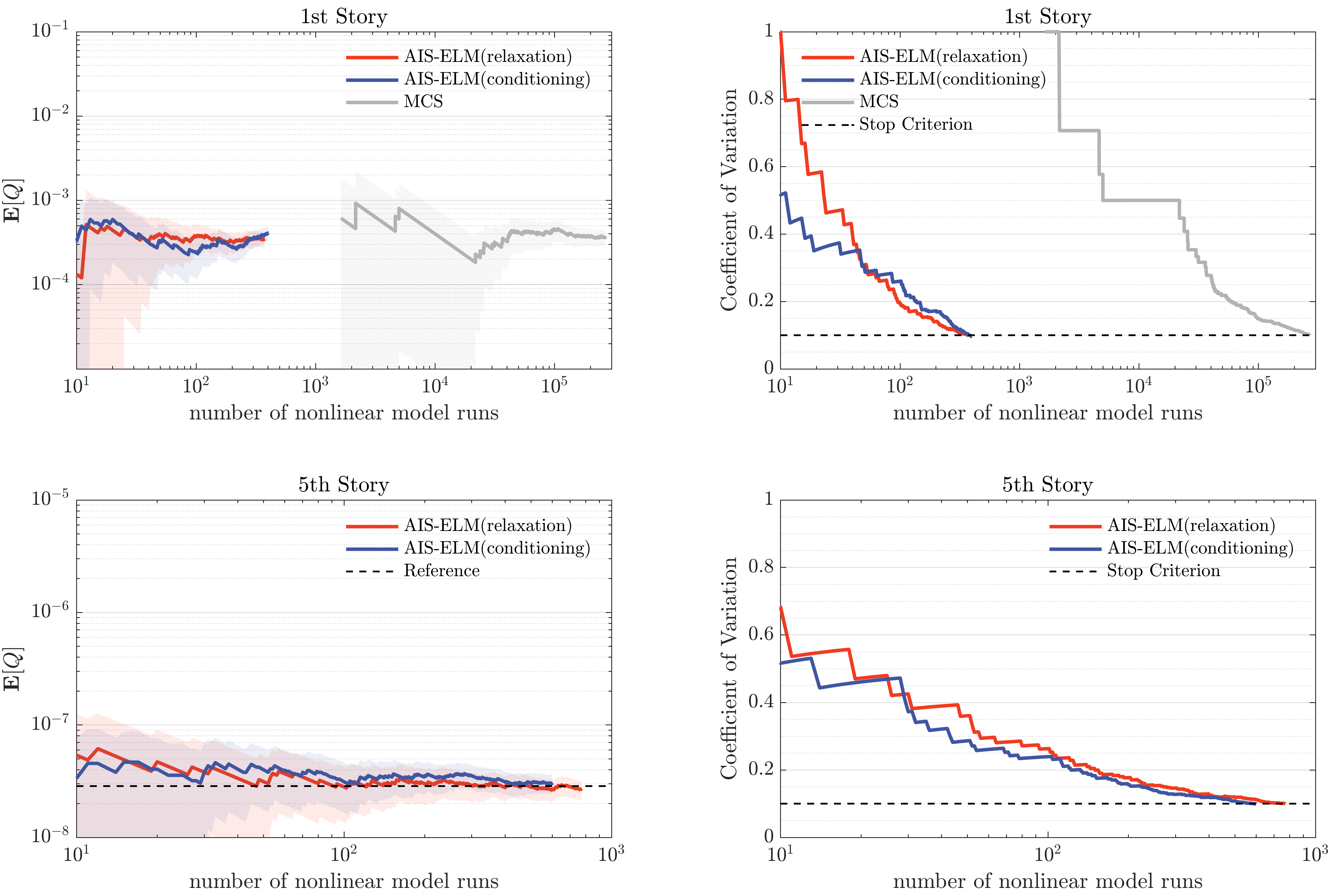}
    \caption{\textbf{Simulation result for the first-passage probability of the 1st and 5th story deformation}. \textit{To achieve the same coefficient of variation of $10\%$, for the 1st story the AIS-ELM based on relaxation requires $686$(for optimization)$+372=1058$ nonlinear model runs, and the probability estimate is $3.48\times10^{-4}$; the AIS-ELM based on conditioning requires $686$(for optimization)$+846=1532$ nonlinear model runs, and the probability estimate is $4.07\times10^{-4}$; the direct Monte Carlo simulation requires $2.66\times10^5$ nonlinear model runs, and the probability estimate is $3.76\times10^{-4}$. For the 5th story, the AIS-ELM based on relaxation requires $600$(for optimization)$+771=1371$ nonlinear model runs, and the probability estimate is $2.71\times10^{-8}$; the AIS-ELM based on conditioning requires $600$(for optimization)$+1000=1600$ nonlinear model runs, and the probability estimate is $3.05\times10^{-8}$; the reference solution obtained from the Hamiltonian Monte Carlo based Subset Simulation \cite{WANG2019} is $2.86\times10^{-8}$. The AIS-ELM based on conditioning involves ``hidden" nonlinear model runs embedded in the MCMC sampler, which cannot be shown in the figure (but are already included in the reported computational cost). The shaded areas represent approximate confidence intervals obtained from the mean plus and minus $1.96$ standard deviations}.}
    \label{Fig:RunFPMDoF}
\end{figure}

Observing the limitation exposed in Figure \ref{Fig:TSMDoF}, one may conjecture that the optimized ELM may not work well in estimating first-passage probabilities, because the first-passage probability estimation is significantly more challenging than estimating the mean peak response. However, note that to efficiently estimate the mean peak response, the optimized linear system should, to some extent, imitate the \textit{global} trend of the nonlinear system response, while estimating the rare event probability of first-passage only requires a \textit{local} fitting over the narrow critical region of the first-passage event.  The first-passage probability estimations of the 1st and 5th story deformation exceeding a threshold of $0.12\,\rm{m}$, with $S_0=0.15\,\rm{m^2/s^3}$, are illustrated in Figure \ref{Fig:RunFPMDoF}. 

Finally, it is worth mentioning that the response threshold for the first-passage event is set to be $12$ times larger than the yield displacement, and the mean peak acceleration excitation is set to be around the standard gravity. This analysis setting pushes the optimized ELM to the limit to test its capability. In engineering applications, there is typically less nonlinearity than what has been studied here.

\section{Additional Remarks, Limitations, and Future Directions}
\subsection{There is no free lunch}
\noindent Most of the so-called ``smart'' probabilistic analysis methods involve a trade-off between efficiency, accuracy, and generality. For example, a method can be both accurate and efficient by implementing some domain/problem-specific knowledge; this makes the method not general. The conventional ELM trades accuracy for efficiency and generality via establishing global linear models. The proposed optimized ELM trades generality for efficiency and accuracy via leveraging local linear models. The unique feature of the optimized ELM is that the solution, in the mean sense, is guaranteed to converge to the actual response quantity of the nonlinear system. The proposed approach is particularly attractive if some critical response quantities of a nonlinear system are of interest, and one looks for an accurate, efficient, and interpretable computational method.

\subsection{The optimized linear system is local}
\noindent The optimized linear system is not a global surrogate for the underlying nonlinear dynamic model of a random vibration problem. Instead, the optimized linear system is a local or end-to-end surrogate model connecting a specific input-output pair. Consequently, the identified linear model varies with input/output specifications. For example, the linear model would change with the intensity of the excitation, making the application to seismic fragility analysis cumbersome. Moreover, to estimate the non-Gaussian response distribution of a nonlinear system, a sequence of optimized linear systems for specified response thresholds needs to be established. This is  similar to the Tail-equivalent linearization method \cite{fujimura2007tail}. Specifically, the first-passage probabilities $\Prob{\sup_{t\in[0,T]}|X(t)|\geq x_i}$ for a discrete threshold set $\lbrace x_i,i=1,2,...\rbrace$ need to be solved; this involves  repeated optimization to find the threshold-dependent linear systems and importance sampling.  Therefore, it is promising to construct efficient computational approaches that account for variations in the input/output specifications, e.g.,  studying inter/extrapolation techniques in the parametric space of linear models.

\subsection{Efficiency gain of the optimized linearization}
\noindent In general, the efficiency gain of the proposed method depends on the ``similarity" between the QoIs of the nonlinear and optimized linear systems. This similarity is measured as the correlation between QoIs in the control variate-based approach and the sampling variance in the importance sampling-based approach. These similarity metrics--Pearson correlation for the control variate and sampling variance for the importance sampling--are from the mathematical constructions of the specific variance reduction methods. However, it is worth reiterating that when the Control Variate-enhanced ELM is applied to estimating the mean peak responses, the Pearson correlation can have compromised accuracy, and the resulting algorithm can be less efficient compared with estimating other response quantities; this motivates the use of Importance Sampling-enhanced ELM. Moreover, the flexibility of parametric linear physical models places a more fundamental constraint on the efficiency gain of the proposed method. For instance, in estimating the mean peak responses of highly nonlinear hysteretic systems (the last two examples), the correlation coefficients can typically achieve $0.8$, with little hope of improving further. We consider this a fundamental limitation of linearization, i.e., the efficiency gain cannot be further enhanced. This motivates the study of more general physics-based surrogate modeling in the upcoming studies. Moreover, this paper focuses on estimating the mean peak responses and first-passage probabilities; future studies are expected to apply the method to other response quantities.

\subsection{Leveraging linear system theories}
\noindent The proposed method optimizes and uses a linear system model to compute a pre-specified response quantity of interest. However, the optimized linear system encodes more information, such as the power spectral density, than merely providing a response quantity value. A promising line of future research is to develop ways to leverage linear system theories fully. 

\subsection{Trend correction}
\noindent A promising direction for future development is to consider a trend function to account for the bias introduced by plasticity or other nonlinear phenomena. To make the approach general, the sigmoid function, widely used in machine learning, could be an attractive candidate for the trend function. 

\subsection{Static problems and equivalent nonlinear systems}
\noindent A natural extension of the proposed method is on static problems. Moreover, optimized nonlinear systems can also be developed within the current framework. Technically, these extensions only involve redefining a parametrized physical system, while all the other components developed in this paper can be inherited. The generalization to parametrized nonlinear systems can overcome the fundamental limitation of linearization and thus further improve the efficiency of QoI estimations. The challenge will be to manage the trade-off between model complexity and accuracy.

\subsection{Physics-informed surrogate modeling}
\noindent The approach developed in this work seems relevant to the increasingly popular field of physics-informed machine learning and surrogate modeling \cite{cuomo2022scientific}\cite{karniadakis2021physics}. However, the difference is that in a typical physics-informed surrogate model, ``physics" is introduced to regularize the parametric space of neural networks so that the training can be effectively pushed toward the physically admissible domain. In this work, we stick with a real physical model, which is physically admissible by construction. However, integrating the proposed method with physics-informed machine learning can be promising.   

\section{Conclusions}
\noindent In this work, an optimized equivalent linearization method (ELM) is developed for nonlinear random vibration analysis. The method tunes a parametric linear physical system such that a response quantity of the linear system can be maximally correlated to the nonlinear system response. Using control variates and importance sampling theories, the control variate-enhanced ELM and importance sampling-enhanced ELM, in conjunction with their optimization equations, are developed. The former approach is suitable for estimating response statistics, while the latter is desirable for rare event probability estimations. The proposed methodology has the convergence property such that as the linear system simulations are guided by a limited number of nonlinear system simulations, the linear system solution can be corrected toward the nonlinear system solution. This property overcomes the long-standing limitation of conventional ELMs. Three nonlinear random vibration examples are studied to test and illustrate the proposed method. It is found that the method is highly promising even when the nonlinearity is strong. To further improve the efficiency and scope of the method, future works can focus on integrating physics-informed machine learning techniques and extensions toward static problems and parametric nonlinear systems. 

\section*{Acknowledgement}
\noindent
I thank Professor Armen Der Kiureghian for commenting and editing several versions of this paper. I thank Dr. Xuefeng He for the support in producing Figure \ref{Fig:UQ}. 

\bibliography{ELM}

\appendix

\section{Implementation Details}\label{Append:Implementation}
\begin{figure}[H]
    \centering
    \includegraphics[scale=0.8]{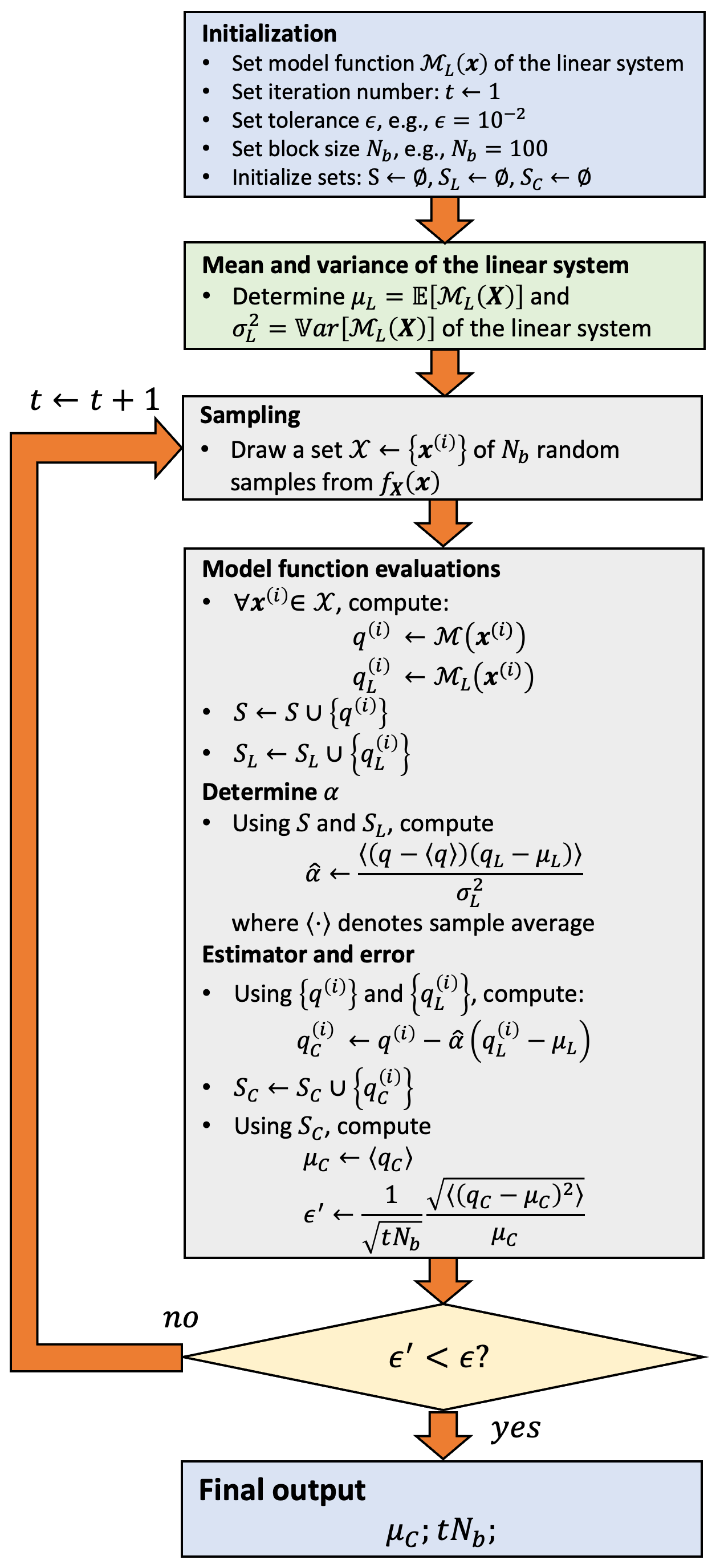}
    \caption{\textbf{Control variate-enhanced ELM}.}
    \label{Fig:CV}
\end{figure}

\begin{figure}[H]
    \centering
    \includegraphics[scale=0.8]{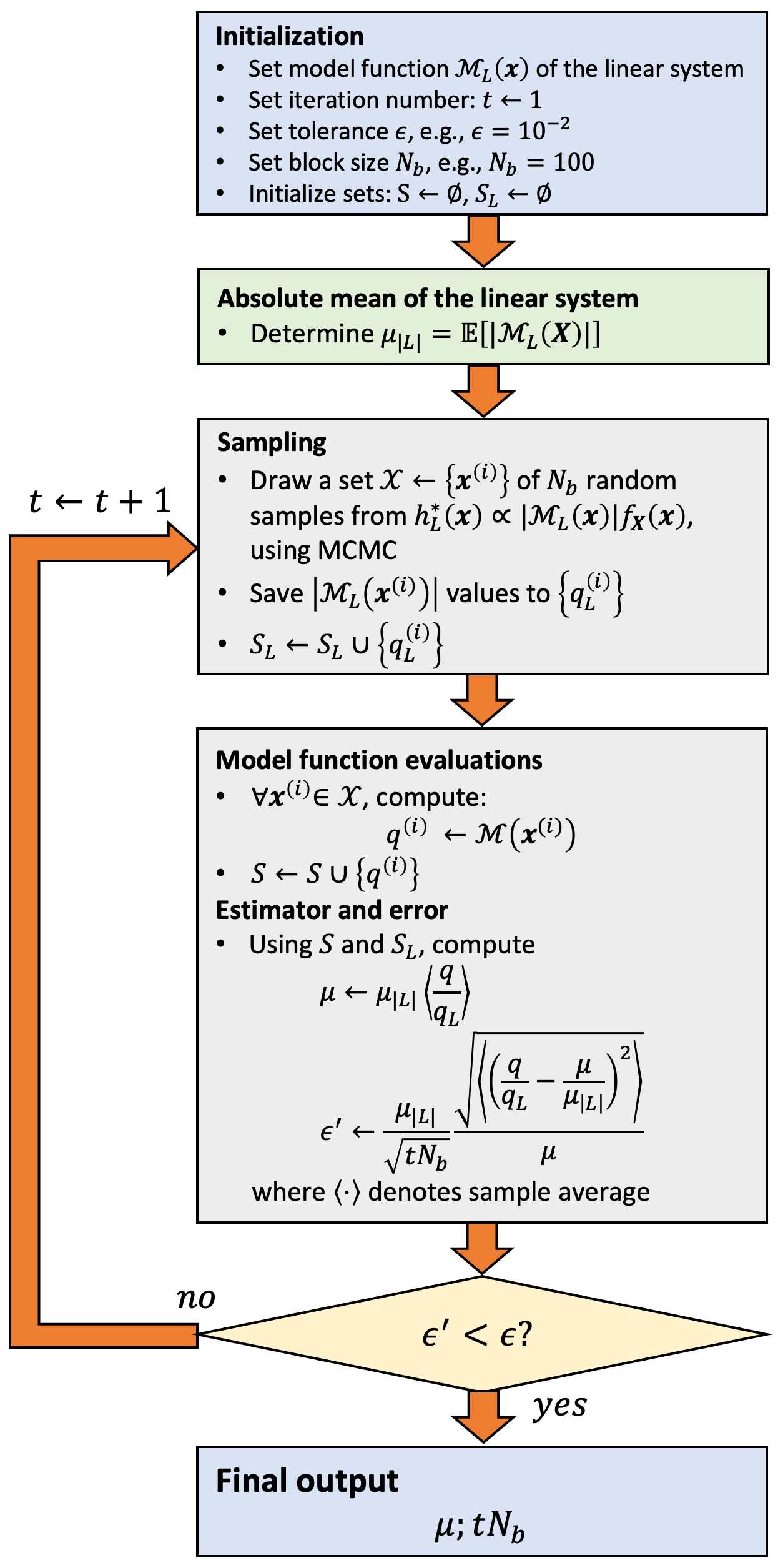}
    \caption{\textbf{Importance sampling-enhanced ELM}.}
    \label{Fig:OIS}
\end{figure}

\begin{figure}[H]
    \centering
    \includegraphics[scale=0.8]{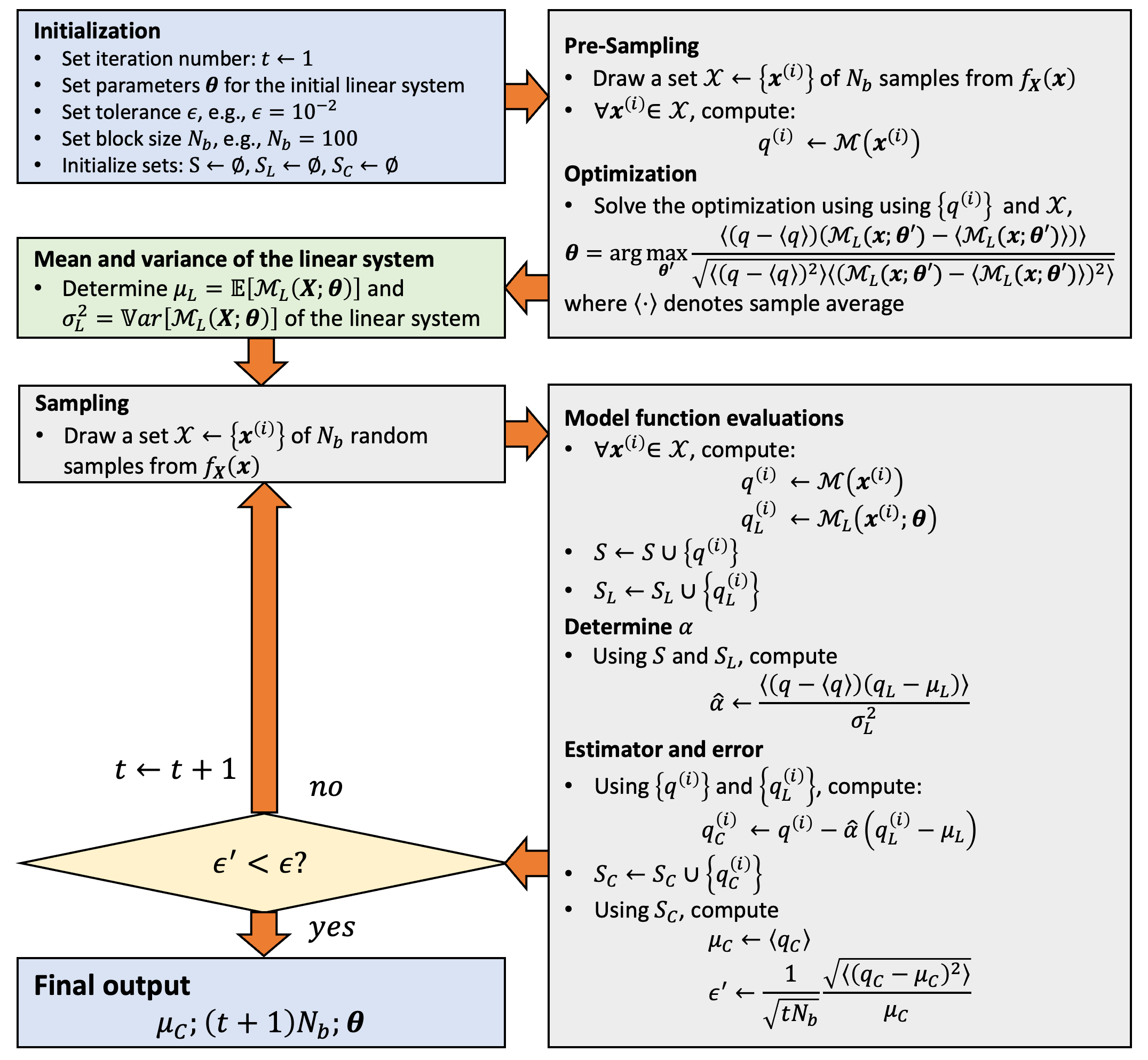}
    \caption{\textbf{Adaptive control variate-enhanced ELM}.}
    \label{Fig:ACV}
\end{figure}

\begin{figure}[H]
    \centering
    \includegraphics[scale=0.8]{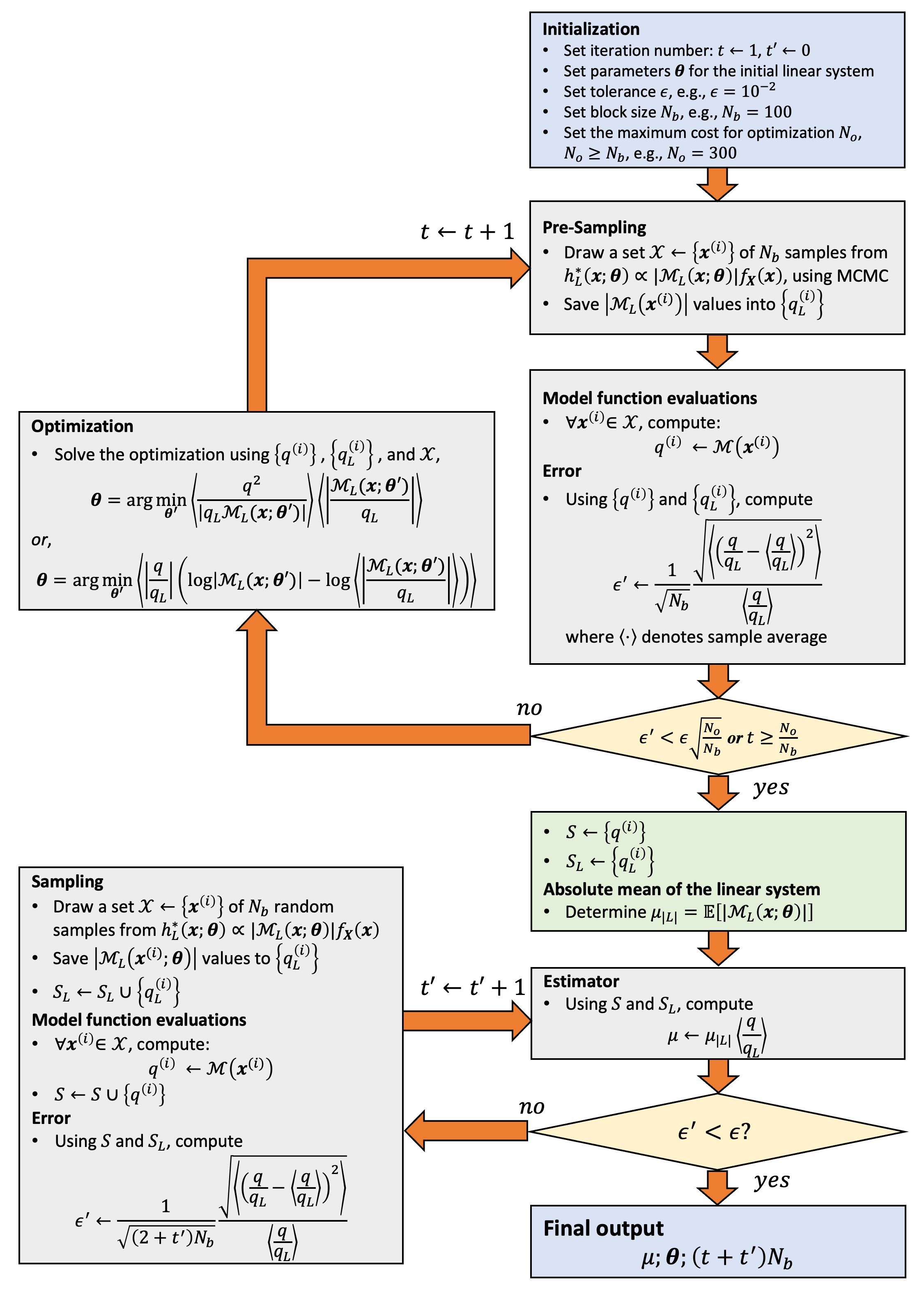}
    \caption{\textbf{Adaptive importance sampling-enhanced ELM}.}
    \label{Fig:AOIS}
\end{figure}

\begin{figure}[H]
    \centering
    \includegraphics[scale=0.8]{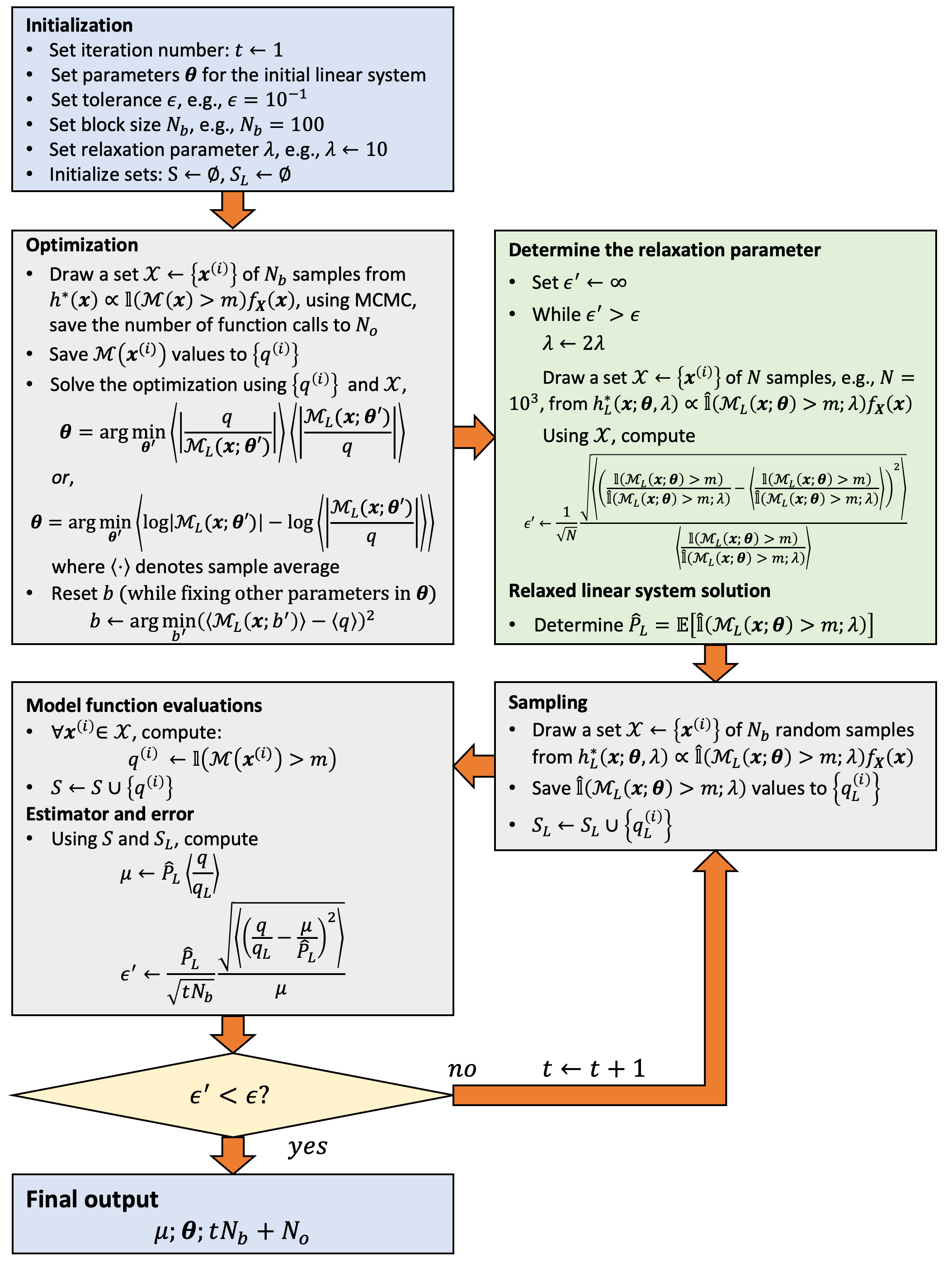}
    \caption{\textbf{Adaptive importance sampling-enhanced ELM (relaxation)}.}
    \label{Fig:ROIS}
\end{figure}

\begin{figure}[H]
    \centering
    \includegraphics[scale=0.8]{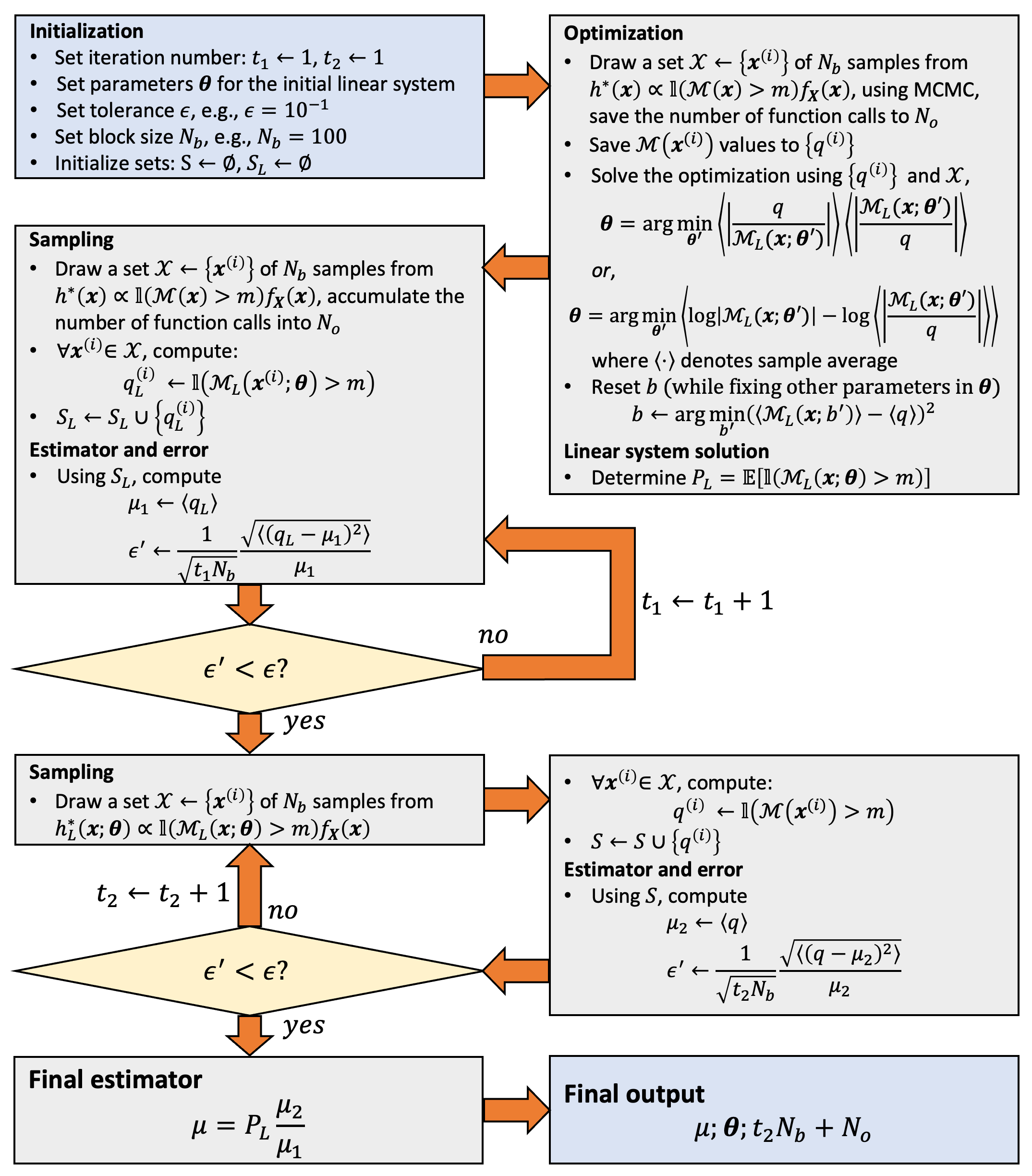}
    \caption{\textbf{Adaptive importance sampling-enhanced ELM (conditioning)}.}
    \label{Fig:COIS}
\end{figure}

\end{document}